\documentclass[12pt]{article}

\input epsf

\usepackage{graphicx}
\usepackage{epstopdf}
\usepackage{amssymb}

\newtheorem{theorem}{Theorem}[section]
\newtheorem{lemma}[theorem]{Lemma}
\newtheorem{proposition}[theorem]{Proposition}
\newtheorem{corollary}[theorem]{Corollary}
\newtheorem{defi}[theorem]{Definition}

\newenvironment{definition}{\begin{defi}\rm}{\end{defi}}

\title{Volume entropy for surface groups via Bowen-Series like maps.}
\author{J\'er\^ome Los, \\
Universit\'e de Provence, LATP,  UMR CNRS 6632.}

\begin{document}

\maketitle

{\bf Abstract.} We define a Bowen-Series like map for every geometric presentation of a cocompact surface group and we prove that the volume entropy of the presentation is the topological entropy of this particular (circle) map. Finally we find the minimal volume entropy among geometric presentations.\\


\section{Introduction}

One proof of the Mostow rigidity theorem for hyperbolic manifolds [Mo] is based on the following variational principle: the manifold admits a unique metric minimising the 
{\em volume entropy} and the optimum is realized by the hyperbolic metric. This approach is due to Besson-Courtois-Gallot (see [BCG] ) and gives some hope for other related area. One such hope would be to obtain an "optimal" metric or an "optimal" presentation in geometric group theory via a similar variational principle. The question of comparing volume entropy in group theory makes sense for Gromov hyperbolic groups (see [Gr1] or [Sho] for definitions), since it is well defined and depends on the presentation. 
The question would be to find a group presentation minimising the volume entropy. An answer to that question is only known for free groups where it is essentially trivial (see for instance [dlH]). 

Recall that a finitely generated group $\Gamma$ with a finite generating set $X$, or a finite presentation 
$P = < X; R> $,  defines a metric space $(\Gamma, d_{X} )$ with the word metric and  $|  B (n) | _{X}$  denotes the cardinality of the ball of radius $n$ centered at the identity in $(\Gamma, d_{X} )$.

The growth properties of the function $ n \mapsto |  B (n) | _{X}$ has been of crucial importance in geometric group theory over the last decades. For instance it is fair to say that a new period started after the fundamental work of M.Gromov classifying groups with polynomial growth functions [Gro2].  Another important step was the discovery by R.Grigorchuck  [Gri] of a whole class of groups with growth function between polynomial and exponential. 

The nature of the growth function (i.e. being polynomial, exponential or intermediate) is a group or geometric property meaning that it does not depends upon the particular presentation (as a quasi-isometric invariant), but the numerical function depends in a highly non trivial way on the group presentation.
For instance, among the exponentially growing groups the following numerical function :\\

\centerline{ $h_{vol} (\Gamma ; P) : = \lim_{n \rightarrow \infty }{ \frac{1}{n} . \ln  |  B (n) | _{X} }$, }

is called the {\em volume entropy} of the presentation $P$ and the way this number varies with $P$ is absolutely not understood.

For hyperbolic groups, the question of finding and characterising minimum volume entropy makes sense, 
but no general method is available to compute or evaluate the volume entropy from the presentation, except for the obvious free group case.

In this paper we develop a method for computing explicitly the volume entropy among the simplest presentations of the simplest (non free) hyperbolic groups, namely among the geometric presentations of co-compact surface groups.

The restriction to geometric presentations is build in our approach. For these special presentations, for which the surface structure is obvious, we re-open a tool box that has been created about 30 years ago by R.Bowen and C.Series [BS]. Their idea was to associate a dynamical system, i.e. 
a $\mathbb{N}$-action on $S^1 = \partial \Gamma $, with all the nicest possible dynamical properties, to a very special group presentation and then to extract some informations about the group from the dynamics.

Bowen and Series defined a Markov map on the circle for one specific presentation of  Fuschian groups, using a special geometric condition on the fundamental domain for the group action in $\mathbb{H}^2$. 
The first new contribution of this paper is to suppress all these geometric conditions but keeping the restriction to geometric presentations. We call a presentation of a surface group 
{\em geometric} if the two dimensional Cayley complex is planar. Generalisation of our construction to arbitrary presentations is conceivable but in a highly non trivial way.
The first result, combining several parts of the paper, can be stated as:

\begin{theorem}
Let $\Gamma$ be a co-compact hyperbolic surface group with a geometric presentation $P$ then there exists a Markov map $\Phi_{P} : \partial \Gamma = S^1 \longrightarrow  \partial \Gamma$ that is orbit equivalent to the group action. In addition this particular map satisfies :\\
\centerline{ Volume Entropy $(\Gamma, P)$  = Topological Entropy  $\Phi_{P}$.}
\end{theorem}

The definition of the map follows the general idea of Bowen an Series but the construction is quite different,  combinatorial here rather than geometric.
The resulting combinatorial dynamical properties of these maps enable to compare the symbolic description of the orbits of the maps with the symbolic descriptions of the geodesics for the given presentation. This comparison is the key step in proving the second part of the Theorem.\\
The Markov map being defined for any geometric presentation, it becomes possible to compare the entropy between different geometric presentations. Furthermore the map is explicit which implies that the computations of the entropy is possible, an exemple is presented in section 5 for the genus two surface. The previous Theorem is the main step in proving the following :

\begin{theorem}
The minimal volume entropy, among all geometric presentations of a co-compact surface group
is realised by the geometric presentations with the minimum number of generators.
\end{theorem}

This result confirms the intuition that presentations with the minimum number of generators are natural candidates to be the absolute minimum for surface groups. The conjecture is still out of reach with the tools developed in this paper. 
Surprisingly, Bowen-Series like maps have been defined in very few cases, Marc Bourdon in his thesis [Bou] and Andre Rocha, also in his thesis [Ro] constructed such maps for some Kleinian groups using a condition for the action of the group on $\mathbb{H}^3$ that is the exact analogue of the condition used by Bowen and Series for Fuschian group actions on $\mathbb{H}^2$.

{\em It is a great pleasure to thank Marc Bourdon and Peter Haissinsky for discussions and comments on this work. }

\section{Some properties of geometric presentations.}

In this section we gather some geometric and combinatorial properties of geometric presentations of surface groups that will be used throughout the paper. 
Recall that a group presentation : $P = \big{ <} x_1, ..., x_n  | R_1, ..., R_k \big{> }$ is given by a set of generators and relations. The relations are words $R_i$ in the alphabet 
$X = \big{\{} x_1^{\pm 1}, ..., x_n^{\pm 1}  \big{\} }$ that are cyclically reduced and are defined modulo cyclic permutations and possibly inversions. The Cayley 2-complex $ Cay^2 (P)$ is the two complex whose 1-skeleton is the Cayley graph $ Cay^1 (P)$ and whose 2-cells are glued to each closed path in the Cayley graph representing a relation. A presentation of a surface group is called {\em geometric} if $ Cay^2 (P)$
is planar. Equivalent definitions that are valid in higher dimension are easy to state (see for instance [FP]).\\
In order to simplify the formulation we assume that the group $\Gamma$ is not a triangular group and has no elements of order two.
In this section most of the statements are easy and are given for completeness.

\begin{lemma}
Let $\Gamma$ be a hyperbolic co-compact surface group and 
$P = \big{ <} x_1, ..., x_n  | R_1, ..., R_k \big{> }$ a geometric presentation of $\Gamma$. Then :\\
1.  The set of generators $\big\{ x_1^{\pm 1}, ..., x_n^{\pm 1} \}$ admits a cyclic ordering that is compatible with the group action.\\
2. There exists a planar fundamental domain $\bigtriangleup_{P}$, where each side $S_i$ of 
$\bigtriangleup_{P}$ is dual to a generator $x_i^{\pm 1}$.\\
3.  Each generator $ x_i $ appears exactly twice (with + or - exponent) on the set of relations 
  $\big\{ÊR_1, ..., R_k  \}$.\\
4.  Each pair of adjacent generators, according to the cyclic ordering (1.), belongs to exactly one relation and defines one relation.

\end{lemma}

These results are classical and can be found, for instance in [FP].$\square$

Let us focus on the boundary of the group.  It is classical that a co-compact surface group, for surfaces of genus larger than 2, is Gromov hyperbolic (see [Gr1] ) and the boundary 
$\partial \Gamma$ is homeomorphic to the circle $S^1 $.
With a presentation $P$, the points $\xi \in \partial \Gamma$
are described as infinite geodesic rays starting at the identity, modulo the equivalence relation, among rays, to be at uniform bounded distance from each other.

These rays are expressed as infinite word representatives, in the alphabet $X = \big\{ x_1^{\pm 1}, ..., x_n^{\pm 1} \}$, considered as infinite paths in the Cayley graph $Cay^1(P)$. We denote $ \{Ê\xi \} $   an infinite word representative of a geodesic ray converging to $\xi \in \partial \Gamma$.
These descriptions are non unique and $\xi \in \partial \Gamma$ have generally more than one geodesic writing. Symbolic description of geodesic rays for surface groups goes back to at least Hedlung in the thirties [He34].

We discuss some properties that are particular to geometric presentations of surface groups.
The non uniqueness of the geodesic writing is reflected by the possible existence of {\em bigons}, i.e. a pair of distinct geodesics $ \{ \gamma_1, \gamma_2 \}$  in $Cay^1(P)$ with the same initial point and the same terminal point. We will often use some classical abuse of notations in identifying the vertices of the complexes $Cay^1(P)$ and $Cay^2(P)$ with the group elements and with some particular writing as geodesic segments ending at those vertices. 

\begin{figure}[htbp]
\centerline{\includegraphics[height=50mm]{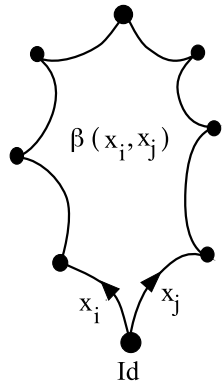} }
\label{fig:1}
\caption{A bigon  in $B(x_i, x_j)$ .}
\end{figure}

Using the group action on the Cayley graph we consider bigons starting at the identity. We denote 
$B(x_i, x_j)$ the set of bigons that start at the identity by the generators $ x_i$ and $ x_j $, for instance 
$\gamma_1 = x_i . w_1$ and 
$\gamma_2 = x_j . w_2$, with $x_i \neq x_j$  (see figure 1).
 
A bigon might be infinite, if the two geodesics $ \{ \gamma_1, \gamma_2 \}$ are geodesic rays, otherwise the length of a bigon is the common length of the two geodesics $ \{ \gamma_1, \gamma_2 \}$. We denote $ \beta (x_i, x_j)$ a bigon in $B(x_i, x_j)$ of minimal length. If necessary we will denote $B_g(x_i, x_j)$ and $ \beta_g (x_i, x_j)$ the bigons based at the vertex $g\in Cay^1(P)$. Observe that the length of a non trivial bigon in a presentation is at least half the minimal length of a relation.

\begin{lemma}
Let $P$ be a geometric presentation of a co-compact surface group $\Gamma$, then :\\
1. $B(x_i, x_j) \neq \emptyset$ only if  $x_i$ and  $x_j$  are two adjacent generators, with respect to the cyclic ordering of Lemma 2.1.\\
2. For each adjacent generators $(x_i, x_j)$ there exists a unique finite length minimal bigon 
$ \beta (x_i, x_j)$.
\end{lemma}
{\em Proof.}
1.  The proof of the first statement is by contradiction. Assume that $B(x_i, x_j) \neq \emptyset$  and 
$( x_i , x_j )$  are not adjacent. If there is a bigon of finite length in  $B(x_i, x_j)$, then we consider a minimal bigon $\beta (x_i, x_j)$. 
The planarity and the minimality assumption imply that $\beta(x_i, x_j)$ is realised by two geodesics whose union is a closed embedded curve in the one skeleton of $Cay^{(2)} (P )$ and therefore bounds a compact topological disc $\mathbb{D}$ in the plane. 
Since $x_i$ and  $x_j$  are not adjacent, according to the planar cyclic ordering of Lemma 2.1, there is at least another generator, say $x'$ between $x_i$ and  $x_j$.  In the Cayley graph, there is a copy of all the generators starting at the vertex denoted $x'$ by an abuse of notation. This vertex and the edges starting at $x'$ are contained in $\mathbb{D}$.
In particular there is another pair of generators $x_i$ and  $x_j$ stating at $x'$ and the geodesics that start by these two edges have to meet, either in the interior of $\mathbb{D}$, or along the boundary i.e. along the paths defining  $\beta (x_i, x_j) $. This intersection defines a bigon in  $B_{x'}(x_i, x_j)$ that is shorter than $\beta(x_i, x_j)$, a contradiction.\\
 For infinite bigons the argument is similar. The two geodesics 
$  \{\gamma_1, \gamma_2\}$ defining the bigon are infinite rays converging towards the same point
$\xi \in \partial \Gamma$. These two rays are disjointly embedded in the plane and bound a disc 
$\mathbb{D}$. They are at distance bounded by some $ \delta $ from each other since the two rays converge to same the point on $\partial \Gamma$. 
The disc $\mathbb{D}$ we consider is such that the restriction of the sphere of radius $N$ with 
$\mathbb{D}$ is a set of diameter bounded by $\delta$.
We assumed that 
$( x_i , x_j )$  are not adjacent hence, there is a copy of the pair $( x_i , x_j )$ at each vertex, at distance one from the identity, beween 
$ x_i $ and $ x_j $ and thus a copy of the disc $\mathbb{D}$ within $\mathbb{D}$. One contradiction comes from the fact that the previous argument implies inductively that the number of vertices on the sphere of radius $N$ within $\mathbb{D}$ grow at least as $3^N$,  a contradiction with the uniform distance between $  \{\gamma_1\}$ and $ \{ \gamma_2\}$.

2.  For the second statement, a pair of adjacent generators $(x_i, x_j)$ defines a unique relation $R$ by Lemma 2.1. This means there is one relation, defined as a cyclic word in the alphabet $X$, that contains the subword $ x_j^{-1} . x_i$ or $ x_i^{-1} . x_j$.

 If the length of the relation $R$ is even :
 
  Then it can be written, up to a cyclic permutation and inversion, as $ w_1.  x_j^{-1} . x_i.w_2 = id$, where the length of $w_1$ and $w_2$ are the same. The two paths written 
$ \gamma_1 = x_j. w_ 1^{-1} $ and $\gamma_2  = x_i . w_2$ connect the identity to the same vertex $z$ in the Cayley graph, where $z$ is the element written as  $ x_j. w_ 1^{-1} $ or $ x_i . w_2$.

{\bf{Claim.}} With the above notations, the two paths $ \gamma_1 $ and $\gamma_2 $ are geodesic segments for the geometric presentation $ P $.

This claim is proved by a contradiction similar to the proof of the first statement. If 
$ \gamma_1 $ and $\gamma_2 $ are not geodesics then there is shorter path $\gamma$ connecting the identity to $z$. The path $\gamma$ has to start with a generator that is different from $x_i$ and $x_j$. The planarity assumption implies that $ \gamma_1 $ lies between $ \gamma_2 $ and $ \gamma $ or 
$ \gamma_2 $ lies between $ \gamma_1 $ and $ \gamma $. In both cases we obtain a contradiction by producing a shorter relation defined by $x_i$ and $x_j$, by the argument of part 1., a contradiction with Lemma 2.1.\\
The pair of geodesics $ \gamma_1 $ and $\gamma_2 $ defines a bigon in $B(x_i, x_j)$ and this bigon is minimal since otherwise there would be another (shorter) relation, defined by the pair $(x_i, x_j)$, a contradiction with Lemma 2.1.

If the length of $R$ is odd:

 The relation can be written as   
$y . w'_1.  x_j^{-1} . x_i.w'_2 = id$, where the length of $w'_1$ and $w'_2$ are the same.
The two paths $ \gamma'_1 = x_j. {w'_ 1}^{-1} $ and $\gamma'_2  = x_i . w'_2$ start at the identity and end at two different points $g_1$ and $g_2$ that differs by the generator $y$. The two paths 
$\gamma'_1$ and $\gamma'_2$ are geodesics by the above argument and $y$ is called "opposite" to the pair  $(x_i, x_j)$.\\
By Lemma 2.1 (item 3.) the generator $y$ appears exactly twice in the set of relations, one of them is
 $R = R^{(1)}$, the other relation $R^{(2)}$ contains the letter $y$ or $y^{-1}$.
If the length of $R^{(2)}$ is odd, then it can be written (modulo cyclic permutation and possibly inversion) as $w''_1. y^{-1} .w''_2 = id$ , where $w''_1$ and $w''_2$ have the same length. The two paths 
$ \gamma''_1 = x_j . {w'_{1}}^{-1}. {w''_{ 1}}^{-1}$ and 
$\gamma''_2  = x_i . w'_2 . w"_2$
are two geodesics from the identity to the same point in $Cay^1 ( P )$ and define a bigon in $B(x_i, x_j)$ that is of minimal length by the above arguments.

If the length of $R^{(2)}$ is even, then it can be written, modulo cyclic permutation and inversion, as  $w'''_1. y^{-1} . w'''_2  . y_1= id$, where $w'''_1$ and $w'''_2$ have the same length. The argument we use for the generator $y$ above is duplicated here for $y_1$. We start an induction on the number of even relations that appear in the following sequence and is uniquely defined from the adjacent pair 
$(x_i, x_j)$:

$ R^{(1)}  \stackrel{\textrm{defines} }{\longrightarrow} {\textrm{opposite : } } \{Êy \}Ê
 \stackrel{\textrm{defines} }{\longrightarrow}  R^{(2)} (\textrm{ even })  \stackrel{\textrm{defines} }{\longrightarrow} {\textrm{opposite : } } \{ y_1  \}\stackrel{\textrm{defines} }{\longrightarrow} R^{(3)} (\textrm{ even }) ....  $\\

If an odd relation appears in the sequence $ R^{(n)}$ then the induction stops because the previous argument defines a unique bigon in $B(x_i, x_j)$ that is minimal for the same reasons.

If no odd relation appear in the sequence then, in particular, $ R^{(1)}$ does not appear again. This implies, in particular, that the generator $y$ does not appear again in the sequence 
$\{ y_{n} \}$. By induction $y_1,...,y_k$ never appear again. This is impossible since the number of generators is finite.  Uniqueness of the minimal bigon is part of the proof. $\square$\\

 \begin{figure}[htbp]
\centerline{\includegraphics[height=35mm,width=150mm]{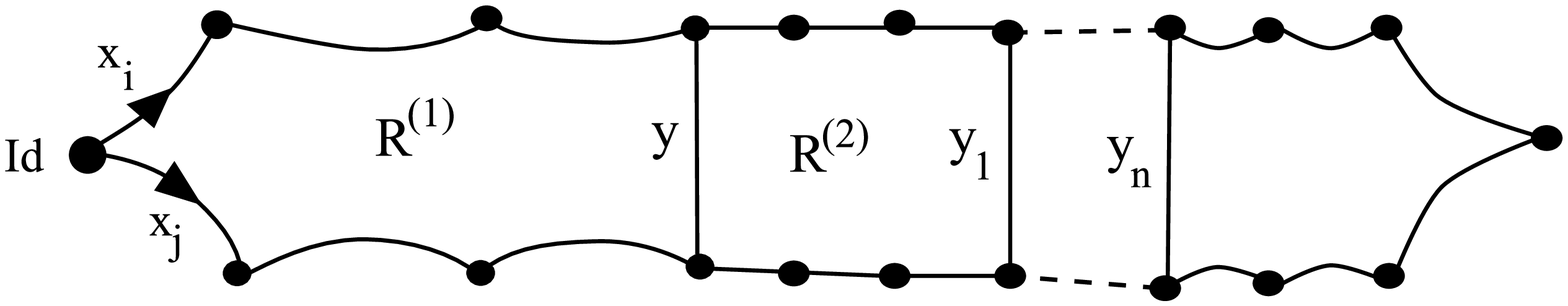} }
\label{fig:2}
\caption{A sequence of relations defining a bigon. }
\end{figure}

A consequence of Lemma 2.2 is :

\begin{corollary}
The boundary $ \partial \Gamma = S^1$ is covered by the cylinders (of length one) 
$C_{x_i}$,\\  $ x_i \in X$, where :\\
$C_{x_i} = \Big\{    \xi \in \partial \Gamma   \textrm{  }  |   \textrm{  } \exists  \{w\} \textrm{  a geodesic ray, representing }  \xi  \textrm{ starting with } { x_i}, \textrm{ i.e. } \{w\} = \{x_i . w'\}   \Big\}$.\\
In addition  $C_{x_i}\bigcap  C_{x_j}  \neq \emptyset $  if and only if $x_i$ and $x_j$ are adjacent generators according to the cyclic ordering of Lemma 2.1.
 \end{corollary}

{\em Proof.} The cylinders cover the boundary since the $x_i$'s generate the group.
A point in $ \partial \Gamma$ belongs to at most two cylinders by Lemma 2.2 (item 1.) and in this case the two cylinders  $C_{x_i}$ and $C_{x_j}$ are defined by two adjacent generators. Conversely the cylinders of two adjacent generators do intersect because of the existence of finite bigons.$\square$

Another consequence of the planarity assumption is:

\begin{lemma}[connectedness]

For a geometric presentation $P$ of a co-compact surface group $\Gamma$, if 
$\xi \in \partial \Gamma $ and $\eta \in \partial \Gamma  $ are two points in the cylinder 
$C_{x_i}$ then one of the two intervals $]  \xi, \eta [ \subset \partial \Gamma $ bounded by $ \xi$ and $ \eta$ is contained in $C_{x_i}$.
\end{lemma}

{\em Proof.} Since $ \xi$ and $ \eta$ belong to $C_{x_i}$ there exists geodesic rays $ \{\xi \}$ and 
$ \{ \eta \}$ starting with $x_i$. The two rays $ \{\xi \}$ and $ \{ \eta \}$ have a common beginning and, since $ \xi$ and $ \eta$ are different, there is a maximal vertex $v$ in $Cay^1(P)$ such that the two infinite paths $ \{\xi \}_v$ and $ \{ \eta \}_v$ starting at $v$ are disjoint. The union $\gamma =  \{\xi \}_v \bigcup \{ \eta \}_v$ is a bi-infinite embedded path in $\mathbb{D}^2$ converging towards $ \xi$ and $ \eta$.
The path $\gamma$ bounds two discs in $\mathbb{D}^2$, one of them contains all the vertices at distance one from the origin, except possibly the vertex corresponding to $x_i$. The other disc, denoted 
$\mathcal{D} (  \{\xi \},   \{ \eta \} ) $, contains one of the two intervals bounded by $ \xi$ and $ \eta$ on the boundary $\partial \Gamma$. We denote $]  \xi, \eta [ $  this interval. Let $\rho \in ]  \xi, \eta [ $, any geodesic ray $\{Ê\rho\} $ representing $\rho$ is contained in $\mathcal{D} (  \{\xi \},   \{ \eta \} ) $ at large distance from the origin. If $\{Ê\rho\} $ has a common initial path with $ \{\xi \}$ or $ \{ \eta \}$ then 
$\rho \in C_{x_i}$. Otherwise, by planarity, $\{Ê\rho\} $ has to intersect either  $ \{\xi \}$ or $ \{ \eta \}$ at a vertex $w$ and therefore defines a bigon in some $B(x_i, x_j)$. In this case $x_j$ is adjacent to $x_i$ by Lemma 2.2 and $\rho$ has another geodesic ray representative $\{Ê\rho\}' $ starting with $x_i$.
$\square$

\section{Special rays and a partition of the boundary.}

In this section we define some special rays in the planar 2-complex $Cay^{(2)} (P)$ giving rise to a finite collection of points on the boundary $\partial \Gamma$ that are uniquely defined from the presentation $P$.
From the previous section, each intersection of two cylinders 
$C_{x_i} \bigcap C_{x_j}$ is non empty only if the two generators $(x_i, x_j)$ are adjacent in $X$. Let us call such a pair of adjacent generators a {\em corner } of the presentation $P$.
The number of corners is even and the cyclic ordering of the generators induces a cyclic ordering of the corners. For notational convenience, the cyclic ordering of the generators is given by the labelling of the generators, in other words $x_{i+1}$ is the generator next to $x_i$ for the cyclic ordering, say on the right.  
By convention it is understood that the notation  $( x_i, x_{i+1} ) $ means that in the 2-complex, the edges denoted $ x_i $ and $ x_{i+1}$ are adjacent and oriented from the vertex. The parity of the number of corners imply that at each vertex, the corner $( x_i, x_{i+1} )$ defines a unique {\em opposite } corner denoted :

(*) {\centerline { $( x_i, x_{i+1})^{opp}\hspace{2mm} :  = \hspace{2mm}(   x_{i +n \hspace{2mm}\mathrm{   mod}[2n]} , x_{i +n  +1 \hspace{2mm}\mathrm{   mod}[2n]} ) $,}}
 
 where $n$ is the number of generators (see figure 3).
 
We construct a unique infinite sequence of corners, bigons and vertices from any given corner
 $( x_i, x_{i+1} )$ by the following process:
 
(i)  Each corner, say at the identity, defines a unique minimal bigon $\beta (x_i, x_{i+1})$ based at $id$, by Lemma  2.2, for which $( x_i, x_{i+1})$ is an extreme corner called the {\em bottom corner}.

(ii)  The bigon $\beta (x_i, x_{i+1})$  has another extreme corner, called a {\em top corner} defined by the end of the two geodesics $\gamma_1 $ and $\gamma_2$ of the definition of a bigon. This extreme corner is denoted:  
$( x_{\beta (i)} , x_{\beta (i+1)} )$ and is based at the vertex $g_1 (x_i, x_{i+1})$. This top corner is uniquely defined by $( x_i, x_{i+1} )$.
 
(iii)  The new corner defines an opposite corner
 $( x_{\beta (i)} , x_{\beta (i+1)})^{opp}$ at $g_1 (x_i, x_{i+1})$.
 
(iv) We consider next the unique minimal bigon: \\
{\centerline {$\beta^{(1)} (x_i, x_{i+1}) := \beta_{g_1} [(x_{\beta (i)} , x_{\beta (i+1)})^{opp}] $,} }\\
that gives a new bottom corner at $g_1 (x_i, x_{i+1})$ and a new top corner at the extreme vertex 
$g_2 (x_i, x_{i+1})$.\\
 This construction defines, by induction, a unique infinite sequence of corners, bigons  and vertices (see figure 3):

\begin{center}
{ $\beta (x_i, x_{i+1})  \longrightarrow    \beta^{(1)} (x_i, x_{i+1}) := 
\beta_{g_1} [( x_{\beta (i)} , x_{\beta (i+1)} )^{opp}]  \longrightarrow    \beta^{(2)} (x_i, x_{i+1})   \longrightarrow \cdots $ }
\end{center}
\begin{center}
 $ ( x_i, x_{i+1} ) \longrightarrow ( x_{\beta (i)} , x_{\beta (i+1)})
 \longrightarrow ( x_{\beta (i)} , x_{\beta (i+1)} )^{ opp} =   ( x_i, x_{i+1} )^{(1)}
 \longrightarrow \cdots   $
 \end{center}

 \begin{figure}[htbp]
\centerline{\includegraphics[height=90mm]{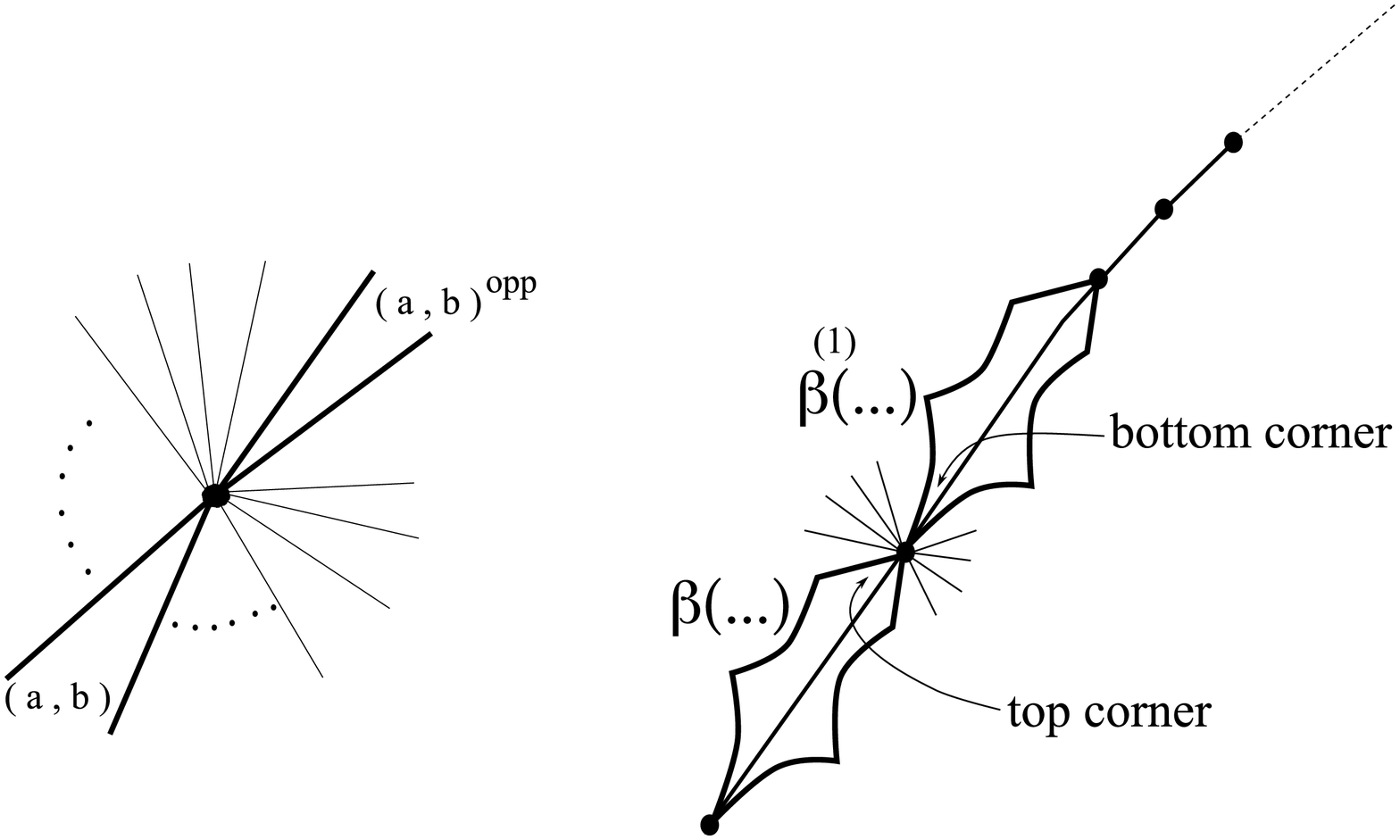} }
\label{fig:3}
\caption{Opposite corner and bigon rays. }
\end{figure}

The bigons that are defined by the previous infinite sequence : $\beta^{(0)} = \beta ,  \beta^{(1)},......,\beta^{(i)}....$, are given by two geodesics $ \big\{  \gamma^{(i)}_{1} ,  \gamma^{(i)}_{2} \}$, $i = 0, 1, 2,...$.

A finite concatenation of bigons :  $\beta^{(0)}.  \beta^{(1)}.......\beta^{(k)}$ is a finite length bigon defined by any finite concatenation of the paths : {\hspace{1cm}
$\gamma^{(0)}_{\epsilon(0)}.  \gamma^{(1)}_{\epsilon(1)}...... \gamma^{(k)}_{\epsilon(k)}$, 
for $ \epsilon(k) = 1 \textrm{ or } 2 $.

\begin{lemma}[bigon rays]
Any of the paths : 
$ \hspace{5mm} \gamma^{(0)}_{\epsilon(0)}.  \gamma^{(1)}_{\epsilon(1)}...... \gamma^{(k)}_{\epsilon(k)}$,  for $ \epsilon(k) = 1 \textrm{ or } 2 $\\
is a geodesic segment in the Cayley graph.\\
In addition, any two such geodesic segments stay at a uniform distance from each other when 
$k \rightarrow \infty$. The infinite concatenation  $\beta^{(0)}.  \beta^{(1)}......\beta^{(i)}.... := \beta^{\infty} (x_i, x_{i+1} )$ is called a {\em bigon ray} and defines a unique point
$( x_i, x_{i+1} ) ^{\infty}$  in  $\partial \Gamma$.
\end{lemma}

The proof of the first statement is a direct consequence of :

\begin{lemma}
Let $P$ be a geometric presentation of a hyperbolic co-compact surface group. If $\gamma$ is a geodesic segment, starting at the identity and ending by a generator $x = x_j \in X$, then any continuation of  
$\gamma$ as  $\gamma . x_i$  is a geodesic segment except when : $x_i = x_j^{-1}$ and possibly when 
$x_i = x_{j \pm1}$ or $x_i = x_{j \pm2}$.
\end{lemma} 

{\em Proof :}  The case $x_i = x_j^{-1}$ is always impossible for a geodesic continuation of the segment 
$\gamma$. The case $x_i = x_{j \pm1}$ might not be a geodesic continuation, in particular at the end vertex of a bigon. The case $x_i = x_{j \pm2}$ might not be a geodesic continuation in the case when the set $R$ of relations contains a relation of length 3. The proof that $\gamma . x_i$ is geodesic in all other cases is obtained by contradiction as in the previous section.
$\square$

The first statement of Lemma 3.1 is obtained inductively using 3.2. Indeed each segment in the sequence is a concatenation of geodesic segments at vertices where the condition of Lemma 3.2 is satisfied by definition (*) of the opposite corner.  Indeed a hyperbolic co-compact surface group has more than 4 generators and, at each vertex of the Cayley graph, more than 8 edges start so the difference of the index by $\pm2$ between the last generator of a segment $\gamma^{(k)}_{\epsilon(k)}$ and the first generator of next segment 
 $\gamma^{(k+1)}_{\epsilon(k+1)}$ is always satisfied.
For the second statement of Lemma 3.1, each minimal bigon in the sequence $\{ \beta^{(n)} \}_n$ has finite length by Lemma 2.2 and the number of different such bigons is finite, this completes the proof of Lemma 3.1.
$\square$

A bigon ray, defined as the infinite concatenation : $\lim_{k\rightarrow \infty} \beta^{(1)} .  \beta^{(2)}... \beta^{(k)}$ is uniquely defined by the corner $ ( x_i, x_{i+1} )$ as well as the limit point on $\partial \Gamma$. This limit point is a particular point in
$ C_{x_i} \bigcap C_{x_{i+1}}$. Each generator $x_i$ belongs to a corner on it's {\em left}
 $( x_{i-1}, x_i )$ and a corner on it's {\em right} $( x_i, x_{i+1} )$ and thus each generator $x_i$ defines two particular points 
$ ( x_{i-1}, x_i )^{(\infty)} $ and $  ( x_{i}, x_{i+ 1} ) ^{(\infty)}$ in $C_{x_i} \subset \partial \Gamma$. Lemma 2.4 implies :

\begin{lemma}
The interval $ I_{x_{i}} : =  [   ( x_{i-1}, x_i )^{(\infty)} , ( x_{i}, x_{i+ 1}) ^{(\infty)}  [  \hspace{2mm} \subset  \partial \Gamma$ is contained in the cylinder $C_{x_i}$ for each $x_i \in X$. $\square$
\end{lemma}

\begin{definition}
For a given geometric  presentation $P$ of a surface group $\Gamma$, the boundary 
$\partial \Gamma = S^1$ admits a canonical partition by the intervals
 $I_{x_{i}}; {\hspace{3mm}} x_i \in X$.
We define the map :

$ \hspace{4cm}$ $ \Phi_{P} : \partial \Gamma \longrightarrow \partial \Gamma$ by 
$\Phi_{P} (\xi ) = x_i^{-1} (\xi)$ when $\xi \in I_{x_{i}}$,

where the action $x_i^{-1} (....)$ is the group action by homeomorphisms on $\partial \Gamma$ induced by the element $x_i^{-1} \in \Gamma$.
\end{definition}

 \begin{figure}[htbp]
\centerline{\includegraphics[height=80mm]{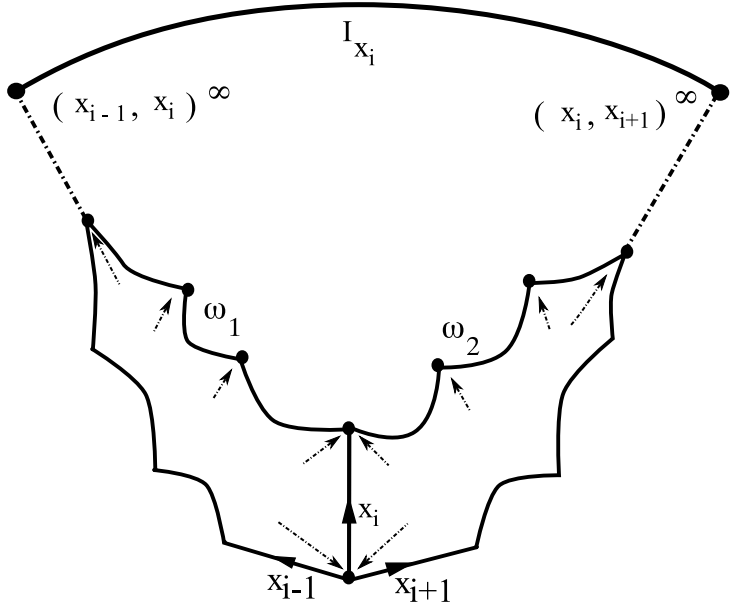} }
\label{fig:4}
\caption{Partition of the circle, definition of $\Phi_{{P}}$. }
\end{figure}

\pagebreak

\section{Subdivision rules and the Markov property.}

The goal of this section is to raffine the partition of $S^1$ by the intervals $I_{x_i}, x_i \in X,$ in order to prove the first part of Theorem 1.1; the Markov property of the map $\Phi_{P}$. 
Recall that a map $F : S^1 \longrightarrow S^1$ satisfies the Markov property if there is a partition (finite here) so that the map is a homeomorphism on each interval and maps extreme points to extreme points. This definition is special to one dimensional spaces (for a more general definition see [Bo] for instance).
From a dynamical system point of view this is just showing that the extreme points of the partition have finite orbits. From a geometric group point of view it is interesting to understand and describe the geometry of the particular geodesic rays that are used to define the partition.

The intervals $I_{x_i}$ are defined through the properties of minimal bigons. The simplest situation is when the presentation has only relations of even length. The next simplest situation is when all the relations are of odd length.
In these simple cases the subdivision process is a little bit easier to describe. 

The interval $I_{x_i}$ is given by the two corners :
$( x_{i-1}, x_{i} )$  and $( x_{i}, x_{i+1} )$.
We focus on the left side of the interval $I_{x_i}$, i.e. on the corner $( x_{i-1}, x_{i} )$, the analysis for the other (right) side is exactly the same. The corner defines a unique relation $R_{L}$, a unique minimal bigon $\beta( x_{i-1}, x_{i})$, a unique bigon ray $\beta^{(\infty)}( x_{i-1}, x_{i} )$
and a unique limit point  $( x_{i-1}, x_{i} ) ^{\infty}$.

\subsection{Simple cases subdivisions.}

We start by assuming that all relations in $P$ have the same parity. The even cases are the simplest situations since all bigons are defined with only one relation (by the proof of Lemma 2.2). In the odd cases all bigons are defined using two relations. All the ideas of the subdivision construction can be seen for these simple cases. \\

(A)  {\bf The length of the relation $R_{L}$ is even:}

In what follows we use a writing that combines some initial path followed by an infinite sequence of bigons :
$\alpha = w . \beta^{\infty}_g (a,b)$, where $w$ is a geodesic path starting at the identity and ending at a vertex $g$ and $ \beta^{\infty}_g (a,b)$ is an infinite sequence of bigons, defined exactly like a bigon ray but starting at $g$ as a geodesic continuation of $w$ and defined by the corner $(a,b)$ at $g$.
This writing describes an infinite collection of geodesic rays. By Lemma 3.1 all the rays in that collection converge to the same point on the boundary.

In particular, among the infinite possible writing, as geodesic rays, of the bigon ray 
$ \beta^{(\infty)}( x_{i-1}, x_{i} )$, the following sub-class makes the belonging of the limit point 
$(  x_{i-1}, x_{i} )^{(\infty)}$ to the cylinder $C_{x_i}$ obvious (see figure 5) :

{\centerline{ $ \beta^{(\infty)}(x_{i-1}, x_{i}  )  \supset  \{  x_{i} . w_{L} .
 \beta^{(\infty)}[(   x_{ \beta (i-1)}, x_{\beta (i)}) ^{(opp) } ]  \}$Ê, where :}} 

$\bullet$  The corner  $( x_{ \beta(i-1)}, x_{\beta(i)} )$ is the top corner of the bigon  
$\beta(  x_{i-1}, x_{i})$.   \\
$\bullet$ The path written $x_{i} . w_{L}=   x_i . y_i^1. y_i^2 .....y_i^k$ is the $x_i$-side of the two paths that define the minimal bigon $\beta(  x_{i-1}, x_{i} )$.

The path $x_{i} . w_{L}$ crosses the following corners :
${ \big\{  { ( \bar{ x_{i}}, y_{i}^{1}),  (\bar{ y}_{i}^{1}, y_{i}^{2} ), ..., (\bar{ y}_{i}^{k-1}, y_{i}^{k} ) \big\} }}$, where $\bar{ z}$ is the standard notation for the inverse orientation.
 Each of those corners   $( a, b)$ defines a unique opposite corner   $( a, b )^{opp}$ and thus a unique bigon ray :  $\beta_{g}^{(\infty)}[( a, b)^{opp}]$ based at the corresponding vertex $g$ (see figure 5).
Starting from the identity we define the following collection of rays, where the based vertex for the bigon rays  $\beta^{(\infty)}$ has been remove to simplify the notations :

[Rays-Even]  {\centerline{$\mathcal{R}_{L} = \Big\{ (  x_i  = y_i^0 ) . y_i^1. y_i^2 .....y_i^j. 
\beta^{(\infty)}  [( \bar{ y}_{i}^{j}, y_{i}^{j+1})^{opp}] ;Ê\hspace{2mm}  j = 0, ... , k-1 \Big\} $. }}

\begin{lemma}
The collection $\mathcal{R}_{L}$ defined above is a collection of geodesic rays called 
{\em left subdivision rays } (with respect to the interval $I_{x_i}$). Each such ray converges toward a point in the interior of $I_{x_i}$.
\end{lemma}

{\em Proof.} The first part is proved the same way than Lemma 3.1.

 \begin{figure}[htbp]
\centerline{\includegraphics[height=90mm]{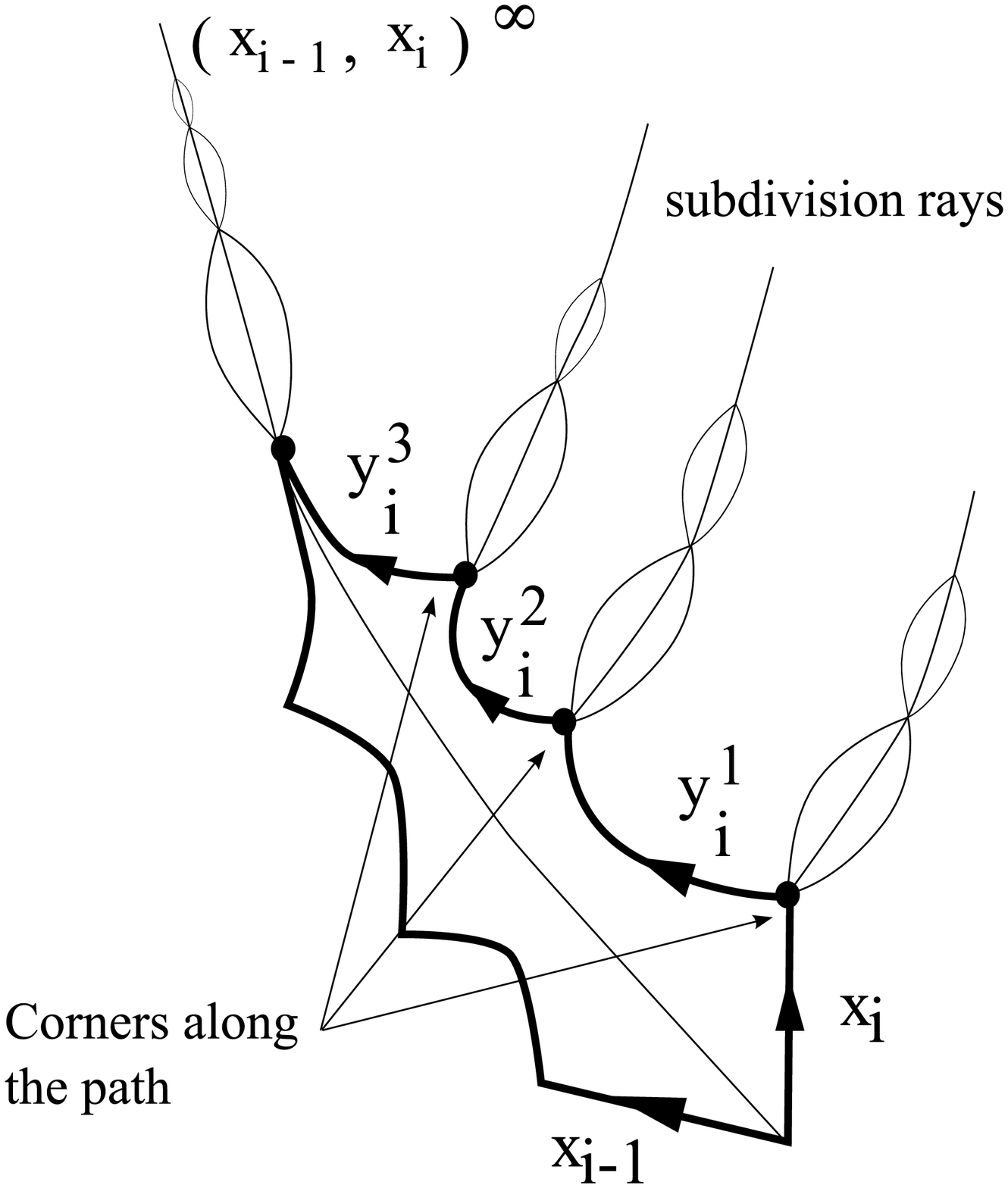} }
\label{fig:5}
\caption{ Subdivision rays, the even case. }
\end{figure}

The collection of rays defining the bigons $\beta^{(\infty)}  ( x_{i-1}, x_{i}  )$ and $ \beta^{(\infty)} ( x_{i}, x_{i+1})$  contain two extreme rays that are disjoint in $\mathbb{D}^2$. The union of these extreme rays is a bi-infinite embedded geodesic, passing through the identity, and bounding a maximal domain $\mathcal{D}_{x_i}$ in 
$\mathbb{D}^2$ with $\mathcal{D}_{x_i} \bigcap  \partial \Gamma  = \overline{ I_{x_i} }$.
 
Any ray that stay in the interior of  $\mathcal{D}_{x_i}$ converges to a point in the interior of  $I_{x_i}$.

By construction the rays in $\mathcal{R}_{L}$ stay in the interior of  $\mathcal{D}_{x_i}$ and thus converge to points in the interior of the interval $I_{x_i}$, that we call {\em (left) subdivision points} 

{\centerline{ ${  \mathbf{ L}}_{x_i} = \{ L_{x_i}^{(1)}, ...,  L_{x_i}^{(k)} \} \subset  I_{x_i} $. 
 $\square$ } }

(B) {\bf The length of $R_{L}$ is odd:}

In this case the definition of the subdivision rays is a little bit more difficult but the idea is just the same.
Recall that in this paragraph all relations are odd. 
The proof of Lemma 2.2 shows that the corner  $ ( x_{i-1}, x_{i} ) $ in the relation $R_{L}$ defines an edge opposite to the corner (denoted $y$ in Figure 2 and $y_i^3$ in Figure 6) on which another relation is based to define the bigon  $\beta( x_{i-1}, x_{i}) $. 

This relation and the corresponding 2-cell in the Cayley complex, is called the {\em bigon completion} of the corner or of the corresponding edge. In the more general case where even and odd relations exist, the single relation is replaced by a unique sequence of relations defining the bigon, the {\em bigon completion} in this general case is this particular sequence.

 The bigon completion is well defined from the corner or from the opposite edge. In the simple case of this paragraph, the bigon completion consists of a single relation.\\
The bigon ray is described, from the $x_i$ side exactly like above (see Figure 6), as :\\

{\centerline{ $\beta^{(\infty)} ( x_{i-1}, x_{i}  ) \supset  \{  x_{i} . w_{L} .   \beta^{(\infty)}[( x_{ \beta (i-1)}, x_{\beta (i)} ) ^{(opp) } ] \} $Ê, where :}} 

$\bullet$  The corner  $( x_{ \beta(i-1)}, x_{\beta(i)} )$ is the top corner of the bigon
 $\beta(  x_{i-1}, x_{i} )$. The top corner belongs to the bigon completion.

$\bullet$  The $x_i$-side of the bigon $\beta( x_{i-1}, x_{i})$ is written as $x_{i} . w_{L} $. This path is expressed as :\\
$   x_{i} . w_{L} =  x_i . y_i^1. y_i^2 .....y_i^m . z_i^{m+1}.... z_i^{m+r}$, where the  first part :   
$x_i . y_i^1. y_i^2 .....y_i^m$  is the $x_i$ side of  $\beta( x_{i-1}, x_{i})$ along the relation $R_{L}$, and the second part :
$ z_i^{m+1}.... z_i^{m+r}$, is the part of the $x_i$ side of the bigon $\beta(  x_{i-1}, x_{i} )$ that belongs to the bigon completion (see Figure 6).\\ 

The path $x_{i} . w_{L} $ crosses the following corners along the relation $R_{L}$ :\\
${ \big\{  { ( \bar{ x_{i}}, y_{i}^{1} ),  (\bar{ y}_{i}^{1}, y_{i}^{2}), ..., 
(\bar{ y}_{i}^{m}, y_{i}^{m+1} ) \big\} }}$. This last corner
 $(\bar{ y}_{i}^{m}, y_{i}^{m+1} )$ is the one that corresponds to the last edge of the path 
 $x_{i} . w_{L} $ along the relation $R_{L}$  ( i.e. ${ y}_{i}^{m}$) and the next one along 
 $R_{L}$ (i.e. ${ y}_{i}^{m+1}$) is the edge opposite to the corner 
 $ ( x_{i-1}, x_{i} )$ for the relation $R_{L}$ (see Figure 6).
Each of those corners $( a, b)$ defines a unique opposite corner
$( a, b)^{opp}$ and thus a unique bigon ray :  $ \beta^{(\infty)}[( a, b)^{opp}] $ based at the corresponding vertex (see Figure 6).
 
\begin{figure}[htbp]
\centerline{\includegraphics[height=90mm]{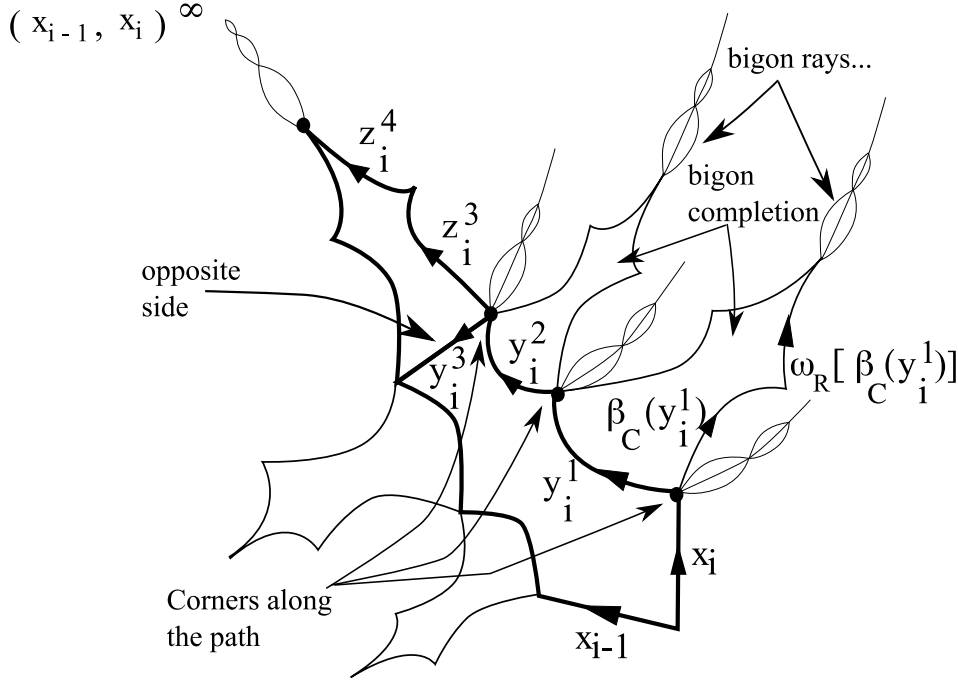} }
\label{fig:6}
\caption{ Subdivision Bigon rays, odd case. }
\end{figure}

The path $   x_{i} . w_{L} =  x_i . y_i^1. y_i^2 .....y_i^m .  z_i^{m+1}.... z_i^{m+r}$ also crosses the edges $\{  y_i^1,  y_i^2,  ..., y_i^m \}$ of the relation $R_{L}$ and are thus the opposite edge of some corner along $R_{L}$.
Based on each of these edges $ y_i^j$ there is a bigon completion of the corresponding opposite corner. We denote $\beta_C (  y_i^j) $ the bigon completion based at the edge $ y_i^j$ and  
$w_{R} [\beta_C (  y_i^j) ]$ the subpath on the right of the bigon completion $\beta_C (  y_i^j) $, 
 up to the top corner of the bigon completion. We denote this corner  $ < \beta_C (  y_i^j)  > $ (see Figure 6), it defines a unique bigon ray  $\beta^{(\infty)} (< \beta_C (  y_i^j)  >^{opp} ) $ based at the corresponding vertex.

Starting from the identity we define the following collection of rays :

$ {\textrm{  [Rays-Odd] }}   \hspace{1.3cm}  \mathcal{R}_{L} = 
\Big\{ (  x_i  = y_i^0 ) . y_i^1. y_i^2 .....y_i^{j}.  \hspace{1mm}  \beta^{(\infty)}[( \bar{ y}_{i}^{j}, y_{i}^{j+1})^{opp}]  ;  \hspace{1mm}  j = 0, ... , m, Ê \Big\}   \bigcup $\\
  ${ \hspace{2.5cm}} $     $ \Big\{  (  x_i  = y_i^0 ) . y_i^1. y_i^2 .....y_i^{j} . w_{R} [ \beta_C (  y_i^{j+1}) ] .    
  \beta^{(\infty)} (< \beta_C (  y_i^{j+1})  >^{opp} ) ;  \hspace{1mm}  j = 0, ... , m-1 \Big\} $. 
 
This notation is not easy to manipulate, we verify, for instance, that the bigon completion  
$ \beta_C (  y_i^{m+1} ) $ of the edge  $ y_i^{m+1}$ is the bigon completion of the original corner 
$( x_{i-1}, x_i ) $. Then the writing :\\ 
{\centerline{ $ x_i  . y_i^1. y_i^2 .....y_i^{m} . w_{R} [ \beta_C (  y_i^{m+1})] . 
   \beta^{(\infty)}(< \beta_C (  y_i^{m+1})  >^{opp} ) $}}\\
  is simply the bigon ray : $  \beta^{(\infty)} ( x_{i-1}, x_i  ) $.
   In particular the path 
  $w_{R}[ \beta (  y_i^{m+1}) ] $ was written above as : $ z_i^{m+1}.... z_i^{m+r}$.

  Just like Lemma 4.1, the collection  $ \mathcal{R}_{L}$  is a collection of rays called ( left)
{\em subdivision rays } (with respect to the interval $I_{x_i}$).
These subdivision rays stay in the domain $\mathcal{D}_{x_i}$ defined above.\\
Therefore all the left subdivision rays converge towards points in the interior of the interval $I_{x_i}$  and are called (left) {\em subdivision points} 
 
    {\centerline {${  \mathbf{ L}}_{x_i} = \{ L_{x_i}^{(1)}, ...,  L_{x_i}^{(2m+1)} \} \subset I_{x_i}$.}}

Observe that in both cases of even and odd relations the set of subdivision points depend only on the combinatorial properties of the relation $R_{L}$.
The description for the right side of  $I_{x_i}$ is the same than for the left side. The only change is to replace the subpath  $ w_{R} [... ]$ on the right by $ w_{L} [... ]$ on the left.
We denote 
${  \mathbf{R}}_{x_i} = \{ R_{x_i}^{(1)}, ...,  R_{x_i}^{(k)} \}$ the set of right subdivision points of 
$ I_{x_i}$.

\begin{lemma} 
The set of subdivision points ${  \mathbf{R}}_{x_i}$ and $ {  \mathbf{ L}}_{x_i}$ of $ I_{x_i}$ are well defined and belong to the interior of the interval $I_{x_i}$. They depend only on the pairs of right and left relations that are uniquely defined by the two corners $( x_{i-1}, x_i ) $ and $( x_{i}, x_{i+1} ) $. 
\end{lemma}

{\em Proof.} The proof is exactly the same than for Lemma 4.1.$\square$

Let us now consider the set of  all subdivision points :
$ \mathcal{S} = \bigcup_{x_i \in X}({\mathbf{R}}_{x_i}\bigcup { \mathbf{ L}}_{x_i} \bigcup {\partial I_{x_i}})$.
 In the simple case of this paragraph, the main technical result is the following:

\begin{proposition}
If the geometric presentation $P$ of the group $\Gamma$ is such that all the relations have the same parity then the partition points $ \mathcal{S}$ are uniquely  defined by the presentation $P$ and are invariant under the map $\Phi_{P}$.
\end{proposition}

{\em Proof.}
We start the proof of 4.3 in the simplest case where all the relations are even. Let $\zeta \in  \mathcal{S}$ be a subdivision point. It is defined by a relation (say on the left )  $R_{L}$.
Such a point is defined by a ray in some 
$\mathcal{R}_{L}$  or is a boundary point in $\partial I_{x_i}$, in both cases and is written as:\\
{\centerline {$\{ \zeta \} =   x_i   . y_i^1. y_i^2 .....y_i^j.  \beta^{\infty}[( \bar{ y}_{i}^{j}, y_{i}^{j+1})^{opp}] $ for some $  j \in \{ 0, ... , k \}  $.} }
Since $\zeta \in I_{x_i}$ then the image under $\Phi_{P}$ is obtained by the action of $x_i^{-1}$ and is represented as the limit point in $\partial \Gamma$ of the rays whose writing are :\\
{\centerline {  \{ $\Phi_{P}(\zeta) \} = y_i^1. y_i^2 .....y_i^j.  \beta^{\infty}[( \bar{ y}_{i}^{j}, y_{i}^{j+1})^{opp}] $ if $j\neq 0$  and }Ê}\\

 {\centerline { $ \{ \Phi_{P}(\zeta) \} = \beta^{\infty}  [ (  ( \bar{ y}_{i}^{0} = x_i^{-1} ), y_{i}^{1})^{opp} ] $ if $j = 0$.} }

The interpretation of this writing is simple. The action of the map, for this particular writing, is a shift map.
Recall that the path $ x_i  . y_i^1. y_i^2 .....y_i^j....y_i^k$ is the $x_i$-side of the bigon 
$\beta ( x_{i-1}, x_{i}  )$ defined by the corner $ ( x_{i-1}, x_{i} )$ in the relation $R_L$.

If $j \neq 0$ then the path $\{ \Phi_{P}(\zeta) \}$ starts by $y_i^1. y_i^2 .....y_i^j$ that begins with  $y_i^1$ at the corner $ ( x_i^{-1}, y_i^1 ) $. The fundamental observation is that the corner 
$ ( x_i^{-1}, y_i^1 ) $ and the path $y_i^1. y_i^2 .....y_i^j$ belong to the same relation $R_L$ (see Figure 7).  In other words the rays in  
$y_i^1. y_i^2 .....y_i^j.  \beta^{\infty}[( \bar{ y}_{i}^{j}, y_{i}^{j+1})^{opp}] $ define one of the subdivision ray, for the interval $I_{y_i^1}$.

If $j =0$ then $ \{ \Phi_{P}(\zeta) \}$ is represented by $   \beta^{\infty} [(  x_i^{-1} , y_{i}^{1})^{opp}]$, this is one of the bigon rays and thus the point $  \Phi_{P}(\zeta) $ is a partition point.
This completes the proof in the cases of even relations.

\begin{figure}[htbp]
\centerline{\includegraphics[height=90mm]{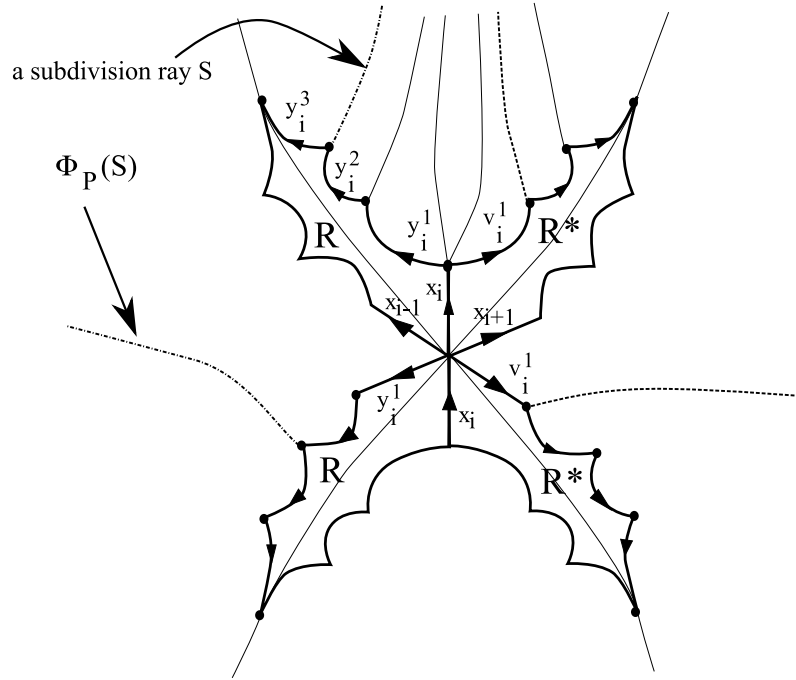} }
\label{fig:7}
\caption{ The action of the map $\Phi_{\Gamma_P}$. }
\end{figure}

In the case where all the relations are odd the above arguments are valid among the collection [Rays-Odd] with one exception for the special ray written above as:\\
$\{  \zeta \} =  x_i . w_{R} [\beta_C (  y_i^{1}) ] . \beta^{\infty} [( \beta_C (  y_i^{1})  )^{opp} ] $
(see Figure 6).
The image of that ray, under $ \Phi_{P}$,  starts with the first letter $w_{R}^1$ of the word 
$ w_{R} [\beta_C (  y_i^{1}) ]$ that is the right side of the bigon completion of the edge $ y_i^{1}$.\\
This bigon completion and the path $ w_{R} [\beta_C (  y_i^{1}) ]$ belong to another relation $\hat{R}$
given by the corner adjacent to $ ( x_i^{-1}, y_i^1) $.
The relations $\hat{R}$ has odd length and the ray :\\
 $w_{R} [\beta_C (  y_i^{1}) ] . \beta^{\infty} (< \beta (  y_i^{1})  >^{opp} )  = \{ \Phi_{P}(\zeta) \} $
 is a subdivision ray defining a subdivision point in the interval $I_{w_{R}^1}$. This completes the proof in the odd case. $\square$

\subsection{Subdivisions for general presentation.}

In order to suppress the parity assumption of the previous paragraph, we need more subdivisions. The reason for this is observed in the last argument of the previous paragraph.

For a presentation with relations of mixed parity, the additional difficulty comes from the minimal bigons with relations of mixed parity. In this case the bigon completion contains some relations of even length and the last argument above produces a ray of the form : \\
$w_{R} [\beta_C (  y_i^{1}) ] . \beta^{\infty} (< \beta_C (  y_i^{1})  >^{opp} )$, where 
$w_{R} [\beta_C (  y_i^{1}) ]$ is a path on one side (right here) of the bigon completion. The path $w_{R} [\beta_C (  y_i^{1}) ]$ starts along an even relation and is followed by a path along some other relations, by the proof of Lemma 2.2.
This ray does not belong to the collection of subdivision rays defined so far.

We need to add more subdivision rays that contain these additional rays.
We consider all the minimal bigons $ \beta(  x_{i-1} , x_i  ) $ of the presentation, this is a well defined finite collection.
The subdivisions defined above at the corners of the odd relations are still given by the collection [Rays-Odd]. We subdivide the same way the intervals with an even relations on one side when it is necessary, i.e. at each edge on which a bigon completion is based. We define the new subdivision rays exactly as in [Rays-Odd] in these cases. We denote again $ \mathcal{S}$ the set of subdivision points of 
$S^1 = \partial \Gamma$ and the same arguments proves the following:

 \begin{proposition}
For any geometric presentation $P$ of the group $\Gamma$, the partition points $ \mathcal{S}$ are uniquely defined by the presentation $P$ and are invariant under the map $\Phi_{P}$. $\square$
\end{proposition} 

\begin{corollary}
The map $\Phi_{P}$ satisfies the Markov property and the collection of subdivision points 
$ \mathcal{S}$ defines a Markov partition for $\Phi_{P}$.
\end{corollary}

The Markov property is clear; on each interval of $S^1 -\mathcal{S}$  the map is the restriction of the action of a group element  and thus is a homeomorphism.  Proposition 4.4  implies that any boundary point of the partition is mapped to a boundary  point of the partition. $\square$

\subsection{Orbit equivalence.}

The orbit equivalence is a relation between the group action on the boundary and the dynamical properties of the map $\Phi_{P}$, more precisely : 

Two points $\zeta$ and $\eta$ in $\partial \Gamma$ are in the same $\Gamma$-orbit if there exists 
$\gamma \in \Gamma$ such that $\zeta = \gamma (\eta)$. The two points are in the same $\Phi_{P}$-orbit if there exists $n, m \in \mathbb{N}$ so that $ \Phi_{P}^n (\zeta) =  \Phi_{P}^m (\eta)$ .\\
The two actions $\Gamma$ and $\Phi_{P}$ on $\partial \Gamma$ are {\em orbit equivalent} if  every points in the same 
$\Gamma$ orbit are in the same $\Phi_{P}$ orbit and conversely.

\begin{theorem}
The two actions $\Gamma$ and $\Phi_{P}$ on $\partial \Gamma$ are  orbit equivalent.
\end{theorem}

One direction of this equivalence is obvious:\\
If $ \Phi_{P}^n (\zeta ) =  \Phi_{P}^m (\eta)$ then the definition of $ \Phi_{P}$ by group elements implies the existence of  $\gamma \in \Gamma$ so that $\zeta = \gamma (\eta)$. 

The other direction requires some work. Assume that $\zeta = \gamma (\eta)$, since $X = \{ x_1^{\pm 1}, ..., x_n^{\pm 1} \}$ is a generating set we restrict to the case  $\gamma =  s \in X $.
Let $\zeta , \eta  \in \partial \Gamma$ be written as limit points of geodesic rays starting at the identity.
 Since $ \zeta = s (\eta )$ then either $\zeta$ admits a geodesic writing as : $ \{\zeta \} = \{ s. x'_{i_2} . x'_{i_3} ... \} $ or $\eta$ admits a geodesic writing $  \{ \eta \} = \{ s^{-1}. x'_{j_2} . x'_{j_3}... \}$. In other words either $\zeta$ belongs to the cylinder $C_s$ or $\eta$ belongs to the cylinder $C_{s^{-1}}$.
Let us assume for instance that $\zeta \in C_s$.

From the definition of $\Phi_P$ then either $\zeta \in I_s$ or $ \zeta \in C_s  - I_s $.
In the first case $\Phi_P (\zeta) = s^{-1} ( \zeta ) = \eta $ and the result is proved. 
In the last case, two situations are possible:\\
(1) $ \zeta \in  I_{x_{i-1}}$, where $x_{i-1} \in X$ is the generator adjacent to $s = x_i$ on the left,\\
(2) $ \zeta \in  I_{x_{i+1}}$, where $x_{i+1} \in X$ is the generator adjacent to $s = x_i$ on the right.\\
We consider only the first case since the two situations are symmetric.
By assumption $ \zeta \in C_{s = x_i} \bigcap C_{x_{i-1}}$ and thus $ \zeta$ admits two geodesic writings :\\
$(*) \hspace{3cm}  \{ \zeta \} = x_i .  w_1  = x_{i-1} . w_2 $.\\
 In addition,  $ \zeta \in I_{x_{i-1}} \bigcap C_{x_{i}}$ and the definition of the subdivision, in section 4.1-2, implies that $ \zeta$ belongs to the right most interval of the partition of $ I_{x_{i-1}}$. This implies that the geodesic writing $(*)$ is given more precisely as :

$(**) \hspace{3cm}  \{ \zeta \} = x_i . w'_1 .  w = x_{i-1} . w'_2.  w  $, where:

$w'_1 = x_{i_1}....x_{i_k} $ and $w'_2  = x_{j_1}....x_{j_k}$ and $\{Ê x_i . w'_1, x_{i-1} . w'_2\}$ are the two geodesic paths defining the minimal bigon $\beta ( x_{i-1} , x_i )$ and $ w $ is a geodesic continuation of $\{Ê x_i . w'_1, x_{i-1} . w'_2\}$ converging to $ \zeta$. This property is a consequence of Lemma 2.2. Indeed, condition $(*)$ implies the existence of a finite length bigon in $B(x_{i-1}, x_i)$ that can be chosen to be the minimal bigon $\beta ( x_{i-1} , x_i )$. \\
By assumption $  \zeta = s (  \eta ) = x_i  ( \eta )$ and thus a geodesic writing of $\eta$ is : \\
{\centerline{ $   \{\eta \} =  w'_1 .  w  = x_{i_1}....x_{i_k} . w  $.}}\\
 The other assumption  $ \zeta \in I_{x_{i-1}}$ implies that : \\
  {\centerline{ $ \{  \Phi_P ( \zeta)  \} =  w'_2.  w  = x_{j_1}....x_{j_k}. w $.}}
  
  We claim that $\eta$ belongs to $I_{ x_{i_1}}$  and $\Phi_P ( \zeta)$ belongs to $I_{ x_{j_1}}$. Indeed the beginning $x_{i_1}....x_{i_k}$ of  $  \{\eta \}$ is a geodesic path starting with $ x_{i_1}$ and is strictly contained in one side of a minimal bigon whereas  $ \{ \Phi_P ( \zeta) \}$ starts with $ x_{j_1}$ by a geodesic path $x_{j_1}....x_{j_k}$ that is strictly contained in one side of another bigon. This implies that the two geodesic rays $  \{\eta \}$ and $ \{ \Phi_P ( \zeta) \}$ belong respectively to the domains
   $\mathcal{D}_{ x_{i_1}} $ and  $\mathcal{D}_{ x_{j_1}} $, as defined in the proof of Lemma 4.1. 
   
   Therefore the $\Phi_P$-image of these two points are :\\ 
{\centerline{ $\hspace{3cm} \{ \Phi_P (\eta) \} =  x_{i_2}....x_{i_k} .   w   $ and 
$ \{\Phi_P ^{2}( \zeta) \}=  x_{j_2}....x_{j_k} .   w   $. }}\\
 By the same argument :  $ \Phi_P (\eta) \in I_{  x_{i_2}}$ and $ \Phi_P ^{2}( \zeta) \in I_{  x_{j_2}}$.  After $k$ iterations we obtain:\\
   {\centerline{ \{ $\Phi_P^{k} (\eta) \} =   w  = \{ \Phi_P ^{k+1}( \zeta) \} $.  $\square$}}

\section{What $\Phi_{P}$ is good for?}

\subsection{Some elementary properties.}

We prove first some simple properties satisfied by $\Phi_{P}$ with non trivial consequences for the group presentation.
Recall that the presentation $P$ defines uniquely the partition $\bigcup_{x_i \in X} I_{x_i}$ and the map 
$\Phi_{P}$ is the piecewise homeomorphism :\\
\centerline{ $\Phi_{P} (\zeta) = x_i^{-1}  (\zeta),   {\hspace{2mm}} \forall \zeta \in I_{x_i } $.}

\begin{lemma}
 If $P$ is a geometric presentation of a co-compact surface group $\Gamma$ with all relations of length greater than 3, let $ I_{x_i}$ be any interval of the partition $ S^1 = \bigcup_{x_i \in X} I_{x_i}$ and
 $ I_{x_{i}^{-1}}$ be the corresponding interval for the inverse generator. If 
  $ I_{({x_{i}^{-1}}) - 1} $ and $ I_{({x_{i}^{-1}}) + 1} $ are the two adjacent intervals  to 
  $ I_{x_{i}^{-1}}$, with $ I_{({x_{i}^{-1}}) - 1} $ on the left (say) and $ I_{({x_{i}^{-1}}) + 1} $
  on the right. Then there exists a subdivision point :
  $\mathbf{ L}_{({x_{i}^{-1}}) + 1}^{j}\in \mathcal{S}$ on the left side of $ I_{({x_{i}^{-1}}) + 1} $
  and another subdivision point:  $\mathbf{ R}_{({x_{i}^{-1}}) - 1}^{k}\in \mathcal{S}$ on the right side of  $ I_{({x_{i}^{-1}}) - 1} $ such that :\\
  {\centerline {$ \Phi_{P}( I_{x_i} )   =  S^1 - [  \mathbf{ R}_{({x_{i}^{-1}}) - 1}^{k} ,  \mathbf{ L}_{({x_{i}^{-1}}) + 1}^{j}  [  $.}}
\end{lemma}
  
 {\em Proof.} Observe that  $ \Phi_{P} $ is a homeomorphism on $ I_{x_i} $ and so it's image is an interval. The Lemma is proved by checking the $ \Phi_{P}$-image of the two extreme points of the interval.   The computation is the same than for the proof of the Markov property (see Figure 7). 
  In particular, in the simple  case where all the relations are even, the subdivision points 
  $\mathbf{ L}_{({x_{i}^{-1}}) + 1}^{j}$ and  $\mathbf{ R}_{({x_{i}^{-1}}) - 1}^{k}$ of the Lemma are  the last point before the $ I_{x_{i}^{-1}}$ and the first point after $ I_{x_{i}^{-1}}$ according to the orientation of the circle. The Lemma is illustrated by Figure 8.
  Let us make the explicit computation in this simple case and for one of the two points. The general case is proved exactly the same way.
   The left boundary point $ (x_{i-1}, x_i )^{\infty}$  of $ I_{x_i} $ is the limit point of the bigon ray   
   $\beta ^{\infty}( x_{i-1}, x_i  )$. In the case where the relation defined by the corner 
   $  (x_{i-1}, x_i )$ is even, this point is written, as in [Rays-Even] : \\
  \centerline{$ ( x_{i-1}, x_i )^{\infty}  =  \Big\{x_i . y_i^1. y_i^2 .....y_i^k. 
\beta^{(\infty)} [( {x}_{\beta(i-1)}, {x}_{\beta(i)})^{opp}]  \Big\} $, }\\
where the path $ x_i . y_i^1. y_i^2 .....y_i^k $ is the right side of the bigon $\beta( x_{i-1}, x_i  ) $.
The definition of $ \Phi_{P} $ on $ I_{x_i} $ implies that the image 
$ \Phi_{P}(  (x_{i-1}, x_i )^{\infty}  )$ is written as the limit point of the ray :\\
\centerline{$  \{ \Phi_{P} [( x_{i-1}, x_i )^{\infty}  ] \}  =    y_i^1. y_i^2 .....y_i^k. 
\beta^{\infty}  [( {x}_{\beta(i-1)}, {x}_{\beta(i)})^{opp}]   $. }

This point is, by definition [Rays-Even], the left most subdivision point of $ I_{ y_i^1} $ (see Figure 7) and this partition interval is adjacent, on the right, to the partition interval  $ I_{x_{i}^{-1}}$. With the notations of Lemma 5.1 we have just checked that:\\ 
 $ \Phi_{P} [(x_{i-1}, x_i )^{\infty}]$ is the last subdivision point, in $\mathcal{S}$, on the left of the interval $ I_{x_{i}^{-1} + 1}$.  The proof for the boundary point on the right is the same and we obtain that   $ \Phi_{P}[(x_{i}, x_{i+1})^{\infty} ]$ is the last subdivision point on the right of the interval $ I_{x_{i}^{-1}- 1}$.  The proof for a general presentations is the same.
 $\square$
  
  \begin{figure}[htbp]
\centerline{\includegraphics[height=60mm]{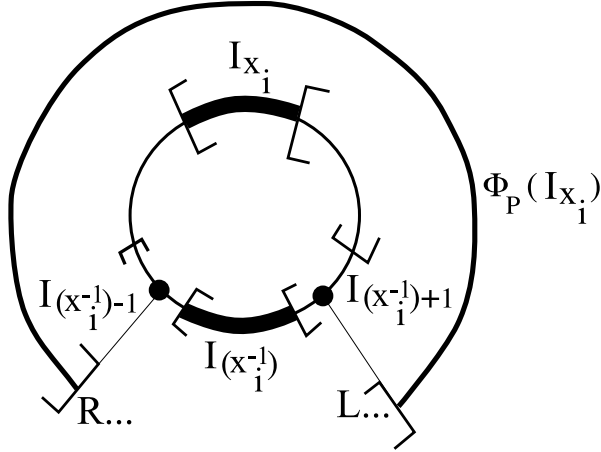} }
\label{fig:8}
\caption{ The image of a partition interval under $ \Phi_{P}$. }
\end{figure}

For presentations with some relations of length 3 the conclusion is a little bit different:

\begin{lemma}
 If $P$ is a geometric presentation of $\Gamma$ with some relations of length 3, for instance $ xyz = Id$, where the 3 generators $\{ x, y, z \}$ are different then, for all generators $x_i$ that do not belong to a relation of length 3, the conclusion of Lemma 5.1 holds.
 For the other generators, for instance $x$ above, then :
 $\Phi_P (I_x ) = S^1 - [ L_a^k , R_b ^j [$, where $ L_a^k$ is a left subdivision point of the interval $I_a$ that is adjacent on the right to $I_{x^{-1}}$ and  $R_b ^j$ is a right subdivision point of the interval $I_b$ that is adjacent to $I_y$  on the left.
\end{lemma}

In other words  the image of the intervals   $I_{x^{\pm 1}}, I_{y^{\pm 1}}, I_{z^{\pm 1}}$ miss two adjacent intervals ( $I_{x^{- 1}}$ and $I_{y}$ in the above case) plus some partition sub-intervals before and after. The proof is the same than for Lemma 5.1 (this particular case is given by Figure 9 ).
$\square$

  \begin{figure}[htbp]
\centerline{\includegraphics[height=90mm]{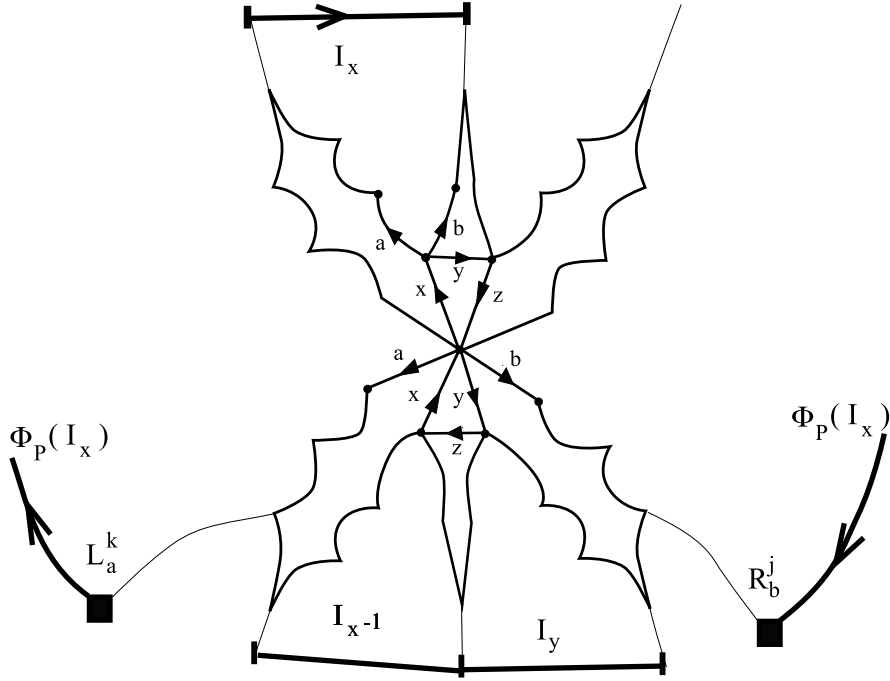} }
\label{fig:9}
\caption{ Image of a partition interval under $ \Phi_{P}$ with a relation of length 3. }
\end{figure}

There is a missing case when the presentation has a relation of length 3 with 2 identical generators. These cases are particular and will be treated in the Appendix. It is interesting to notice that this particular type of presentations exist as a geometric presentation of a surface group but only for non orientable surfaces.

  The two previous Lemmas have the following consequence :
  
  \begin{corollary}
  If $P$ is a geometric presentation of $\Gamma$ with $n$ generators and no relations of the form 
  $xxy = Id$ then 
   the map $ \Phi_{P}$ is strictly expanding and, for all $x\in S^1$, the number of pre-images under
    $ \Phi_{P}$ satisfies :\\
  {\centerline{$ 2n - 3 \leqslant |  \{ \Phi_{P}^{-1} (x) \}Ê |  \leqslant  2 n - 1$.}} 
  
  \end{corollary}
  
{\em Proof.}  The expansivity property is proved directly, each interval of the Markov partition is mapped either on a single interval that is different from itself or is mapped to a union of intervals. In any case the second iterate is mapped to a union of more than two intervals. This argument also proves the transitivity of the map  $ \Phi_{P}$.\\
In the case where $P$ has no relations of length 3 then each point in $S^1$ belongs to one of the intervals $I_{x_i}$ and each such interval has a right and a left subinterval given by Lemma 5.1.
The points in either the right or the left subintervals have exactly $2n - 2$ pre-images and the points in the "central" part have exactly $2n - 1$ pre-images. In this case only the image of the interval corresponding to the inverse generator is not a pre-image. If $P$ has some relations of length 3 with 3 different generators then Lemma 5.2 implies that some points might have $2n - 3$ preimages.
$\square$

Observe that the cases of presentations with relations of length 2 has not been considered. In this case there is a redundant generator that can be removed. This fact is obvious from a combinatorial group theory point of view, it will be justified from a dynamical system point of view in the appendix (Lemma 6.2).

\subsection{Symbolic coding from Markov partition.}

The map $\Phi_{P}$ is Markov and strictly expanding then the standard methods in symbolic dynamics apply (see for instance [Shu], [Bo]) and define a symbolic coding of the orbits :\\
$\{  \Phi_{P}^n (\zeta) \hspace{2mm}; n \in \mathbb{N},  \zeta\in S^1 \}$. As usual in this context we don't consider the finite collection of orbits of the boundary of the partition $\mathcal{S}$.

\begin{definition}
Each point $\zeta \in  S^1 - \mathcal{S}$ admits a symbolic coding on the alphabet 
$ \mathcal{I} =\{Ê I_{(i,j)} \}$, where each $ I_{(i,j)}$ is an interval of the Markov partition of $  \Phi_{P} $, with:\\  
\centerline{ $ \bigcup_{j = 1, ...., K_i} I_{(i,j)} = I_{x_i} \hspace{2mm}, x_i \in X$.}
This coding $ \mathcal{C} :  S^1 - \mathcal{S}  \rightarrow \mathcal{I}^*$, where $\mathcal{I}^*$ is the set of infinite words in the alphabet $ \mathcal{I}$, is defined by :\\
\centerline{$ \mathcal{C}(\zeta) = \{   I_{(i_1,j_1)}, I_{(i_2,j_2)}, ..., I_{(i_k,j_k)} , .... \} \in \mathcal{I}^* $,
where $ \Phi_{P}^k (\zeta)  \in  I_{(i_k,j_k)}, \forall k \in \mathbb{N}$.}
\end{definition}

\begin{lemma}
The coding $ \mathcal{C}$ is called an $I$-coding, it is injective and defines a sub-shift of finite type.
\end{lemma}

The sub-shift property is classical for expanding Markov maps (see for instance [Bo]) as well as the injectivity.  It is a consequence of the fact :\\ 
{\centerline{ ${\textrm{length  }} [ \bigcap_{k=1}^{\infty} \Phi_{P}^{- k}(I_{(i , j)} )]  = 0, {\textrm{for all  }}  I_{(i , j)} $.  $\square$ }}

What is less standard is the fact that the map 
$ \Phi_{P}$ is defined via a group action. In addition, the partition and the action, reflect the action of the generators of the group presentation. The $I$-coding induces an $X$-coding by forgetting the second index of each letter $I(i,j)$ in the alphabet  $ \mathcal{I}$, more precisely:

\begin{definition}
The $X$-coding  $ \chi$ of any $\zeta \in  S^1 - \bigcup_{x_i \in X} \partial I_{x_i}$ is defined  by :\\
\centerline{ $ \chi (\zeta) = \{  x_{i_1}, x_{i_2}, ...., x_{i_k}, ... \}$,
where $ \Phi_{P}^k (\zeta)  \in  I_{x_{i_k}} , \forall k \in \mathbb{N}$.}
\end{definition}

\begin{lemma}
The  $X$-coding is injective and defines a sub-shift of finite type.
\end{lemma}

The $X$-coding is injective, again by expansivity of the map, and the fact that :\\
${\textrm{length  }} [ \bigcap_{k=1}^{\infty} \Phi_{P}^{- k}( I_{x_{i}} )]  = 0, {\textrm{for all  }}  I_{x_{i}} $.
$\square$

For the two codings $ I$ and $X$, any initial word of length $k$ is called a {\em  $I$-prefix (resp. $X$-prefix) of length $k$}.

\subsection{Comparison of entropies.}

The last statement of the first main result, Theorem.1.1,  is a comparison between the asymptotic geometry of the presentation $P$ and the asymptotic dynamics of the map $ \Phi_{P}$.

\begin{theorem}
The volume entropy $h_{\textrm{vol}} (P) $ of the presentation $P$ is equal to the topological entropy 
$h_{\textrm{top}} (\Phi_{P} ) $ of the Markov map  $\Phi_{P}$.
\end{theorem}

 The volume entropy is defined in the introduction and is now a classical invariant in geometric group theory (see for instance [DlH]). The topological entropy of a map is an even more classical topological invariant that was first defined in [AKM] and precised latter by Bowen ( see for instance in [Bo]).
 An important feature of Markov maps is the well known fact that it's topological entropy is computable as the largest eigenvalue of an integer matrix (see bellow).  

One way to relate the geometry of $P$ with the dynamics of $\Phi_{P}$ is to introduce a decomposition of the complexes $Cay^{(j)}(P), {\textrm{for  }} j = 0, 1, 2$  `` suited '' with the intervals $I_{x_{i}}$ on the boundary. This decomposition is a weak combinatorial version of the `` half spaces" in hyperbolic geometry.

Recall that the partition interval $I_{x_{i}} = [ L_i = ( x_{i-1},  x_{i} )^{\infty} ,  R_i = ( x_{i},  x_{i+1} )^{\infty}  [ \hspace{2mm} \in \partial \Gamma $ is defined by two points that are limit of infinite geodesic rays, as in section 3. The bigon ray 
$\beta ^{\infty}( x_{i-1},  x_{i})$ whose limit point is $L_i$ is the concatenation:
  $ \beta_{i,L}^1 .  \beta_{i,L}^2 .....  \beta_{i,L}^k ...$, where the bigon $ \beta_{i,L}^k$ is defined by two geodesic paths : $  \{  ( \gamma_{i,L}^k)^L ,   ( \gamma_{i,L}^k)^R\}$, where the indices $L$ and 
  $R$ stands for left and right. A similar writing defines the point $R_i$.
We define now a particular representative of the geodesic rays converging to $L_i$ and $R_i$ as :\\
\centerline{ $ \{ L_i \}^R : = ( \gamma_{i,L}^1)^R .  ( \gamma_{i,L}^2)^R.....( \gamma_{i,L}^k)^R...$}\\
\centerline{ $ \{ R_i \}^L : = ( \gamma_{i,R}^1)^L .  ( \gamma_{i,R}^2)^L.....( \gamma_{i,R}^k)^L ...$}.

These two particular rays start at the identity by the same letter $x_i$ and are otherwise disjoint.
The compactification $\overline{ Cay^2 (P) }$ of the 2-complex $Cay^2 (P)$ is homeomorphic with the disc $\mathbb{D}^2$ and the union of the two rays $ \{ L_i \}^R \bigcup  \{ R_i \}^L$ is a bi-infinite geodesic connecting $L_i$ and $R_i$, after removing the initial common segment given by the letter $x_i$. This bi-infinite geodesic is an embedded Jordan curve in $\mathbb{D}^2$ and therefore it bounds a domain :\\
$\overline{K_i^2 (P) }$ in $\overline{ Cay^2 (P) }$ so that 
$\overline{K_i^2 (P) }\bigcap \partial \Gamma =  \overline{I_{x_i} }$. This domain is contained in the domain $\mathcal{D}_{x_i}$ of section 4.1.

If we consider the " adjacent'' domains $\overline{K_{i-1}^2 (P) }$ and $\overline{K_{i+1}^2 (P) }$ for the adjacent generators $x_{i-1}$ and $x_{i+1}$ we obtain the following property :

\begin{lemma}
The domains $\overline{K_i^2 (P) }$ defined above satisfy :\\
$\overline{K_i^2 (P) } \hspace{2mm} \bigcap \hspace{2mm} \partial \Gamma = \overline{I_{x_i}} \hspace{5mm}$ and 
$  \hspace{5mm} \overline{K_i^2 (P) } \hspace{2mm}\bigcap  \hspace{2mm} \overline{K_{i-1} ^2 (P) } = \hspace{2mm} L_i \hspace{2mm} \bigcup \hspace{2mm} \{g_{i}^1 ;  g_{i}^2 ; .... \}$,  where \\
$ g_{i}^k$ are the special vertices in $ Cay^2 (P)$ where two minimal bigon paths meet along the bigon ray $\beta ^{\infty}( x_{i-1},  x_{i})$ .
The domain $\overline{K_i^2 (P) }$  is said to be {\em suited } with $I_{x_i}$.
\end{lemma}

{\em Proof.} The bigon ray $\beta ^{\infty}( x_{i-1},  x_{i})$ is the union of the two special rays 
 $ \{ L_i \}^R$ and  $ \{ R_{i-1} \}^L$, where $ \{ L_i \}^R$ is the boundary (left) of $ \overline{K_i^2 (P)}$ and $ \{ R_{i-1} \}^L$ is the (right) boundary of  $ \overline{K_{i-1}^2 (P)}$. These two rays meet at the extreme vertices of the bigons  $ \beta_{i,L}^k$, which are precisely the set 
 $\{g_{i}^1 ;  g_{i}^2 ; .... \}$.$\square$

The set of domains $\overline{K_i^2 (P) }$ for all $x_i \in X$ is not a partition of 
$\overline{ Cay^2 (P) }$ because infinitely many 2-cells are missing, as well as possibly infinitely many 1-cells. For 0-cells, identified with the group elements, the situation is simpler, as given by :

\begin{lemma}
For every $m \in \mathbb{N}^*$, each group element $g \in \Gamma - Id$ of length $m $, with respect to the presentation $P$, belongs to exactly one domain $\overline{K_i^2 (P) }$, except at most $2n$ elements, where $n$ is the number of generators of $P$.

\end{lemma}

{\em Proof.} By planarity, each $g \in \Gamma$ belongs to at most two $\overline{K_i^2 (P) }$ 
that are adjacent. Lemma 5.9 implies that if a vertex belongs to more than one domain it has to be on the boundary of the domain and there are at most two such points on the intersection of the sphere of radius 
$m$ with a domain.
$\square$

From the definition of the domains, each $ z \in \overline{K_i^2 (P) }$ admits a geodesic writing as
$ z  = x_i .  w  $, where the path written $ Êw $ is contained in $\overline{K_i^2 (P) }$ and is possibly infinite.
If $ z = g$ is a group element, i.e. of finite length, there is an open connected set $\Omega_g \subset I_{x_i}$ such that any $\zeta \in \Omega_g$
belongs to the "shadow" of $g$ in  $I_{x_i}$, i.e. $\zeta $ has a geodesic ray writing as 
$\{ \zeta \} =  g .  \rho $  and the geodesic ray is contained in $ \overline{K_i^2 (P) }$. This property comes from the definition of the domains.

For every point $\zeta \in \Omega_g \subset I_{x_i}$ the map $\Phi_P$ is well defined and  
$\Phi_P (\zeta) $ is given by a geodesic ray as 
$\{ \Phi_P (\zeta) \} =  w . \rho \subset  \overline{K_{i_1}^2 (P) } $.
We iterate the argument and there exists a geodesic ray:  
$\{ \Phi_P (\zeta) \} =  x_{i_1}.w ' . \rho $, where
 $ x_{i_1}.w '   \in   \overline{K_{i_1}^2 (P) }$ and $g$ admits a geodesic writing\\
$\{ g \} = x_i .  x_{i_1} .w '  $. After a finite iteration of $\Phi_P $ on $\zeta $ we obtain a geodesic writing of $g$ as $g  = x_i .  x_{i_1}.... x_{i_m}$  and the word : 
$x_i .  x_{i_1}.... x_{i_m}$ is a $X$-prefix of length $m+1$. This proves the following :

\begin{lemma}
Every $g \in \Gamma$ admits a $X$-prefix as a geodesic writing in the presentation $P$.  $\square$
\end{lemma}

The coding by $X$-prefix is not bijective for the group elements  by Lemma 5.9, but Lemma 5.10 gives a uniform bound on the number of elements with more than one $X$-prefix.

We make now a counting argument, let us denote :\\
$\bullet$  $ \sigma_m$ the number of elements of $\Gamma$ of length $m$ with respect to $P$.\\
$\bullet$  $ X_m$ the number of $X$-prefix of length $m$.\\
$\bullet$  $ I_m$ the number of  $I$-prefix of length $m$.\\

\begin{lemma}
For $m\in \mathbb{N}$ large enough, one has $ \sigma_m \approx X_m$.
\end{lemma}

{\em Proof.} If the coding by $X$-prefix were bijective we would have $ \sigma_m = X_m$ for all $m$. This is not the case but Lemma 5.10 implies that these two numbers are equivalent for large $m$. $\square$\\
The last counting argument relates $X$-coding with $I$-coding :

\begin{lemma}
There is a constant $K$ so that for all $m \in \mathbb{N}$ :
$ X_m \leqslant I_m \leqslant X_{m+K}$.
\end{lemma}

{\em Proof.} Every $I$-coding defines an $X$-coding by the  projection :\\
\centerline{ $ I$- coding $ \longrightarrow$  $ X$- coding, given by : }\\
\centerline{ $ \big\{   I_{(i_1, j_1)}, ...., I_{(i_m, j_m)}, ...\big\}  \longrightarrow  \big\{   x_{i_1}, ...., x_{i_m}, ...\big\} $.}
This map is surjective on prefixes by Lemma 5.11, which proves the first inequality.
 Let $K$ be the maximal number of subintervals of the Markov partition among the intervals  $I_{x_{i}}, x_i \in X $. 
If two points $\zeta$ and $\rho$ have the same $I$-prefix of length $m$ then 
$\Phi_P^k (\zeta)$ and $\Phi_P^k (\rho)$ belong to the same $I_{(i,j)}$ for each 
$k \in \{0, 1, ..., m\}$
and in particular they belong to the same $I_{x_i}$. The worst case is when 
$\Phi_P^m (\zeta)$ and $\Phi_P^m (\rho)$ belong to the same $I_{x_i}$ but different $I_{(i,j)}$. By expansivity of $\Phi_P$, after at most $K$ iterations $\Phi_P^{m+K} (\zeta)$ and $\Phi_P^{m+K} (\rho)$ belong to different $I_{x_i}$. So different $I$-prefixes of length $m$ implies different 
$X$-prefixes of length at most $m+K$. This implies the second inequality.
$\square$

{\em Proof of Theorem 5.8.}\\
By definition of the volume entropy of  $P$:\\
$h_{\textrm{vol}}(P) = \lim_{m \rightarrow \infty} \frac{1}{m} \log (\sigma_m)$ and Lemma 5.12 implies  :
 $ h_{\textrm{vol}}(P) =  \lim_{m \rightarrow \infty} \frac{1}{m} \log (X_m)$.\\
  Finally Lemma 5.13 implies :
$h_{\textrm{vol}}(P) =  \lim_{m \rightarrow \infty} \frac{1}{m} \log (I_m)$.

This last limit is classical (see [Shu] for instance) to compute via the "Markov transition matrix" $M(\Phi_P)$ defined as the integer matrix whose entries $M_{a,b} = 1 $ if the interval $I_b$ of the Markov partition (denoted $I_{(i,j)}$ above) is contained in the image $\Phi_P (I_a)$ and $M_{a,b} = 0 $ otherwise.
It is also classical that the norm $\| M ^m \| = \sum_{a, b} | M^m_{a,b} |$ is exactly the number $I_m$ and :

\centerline { $ \lim_{m \rightarrow \infty} \frac{1}{m} \log (I_m) =  \lim_{m \rightarrow \infty} \frac{1}{m} \log (\| M ^m \|)  = \log (\lambda( M) )$, }
where $\lambda ( M)$ is the largest eigenvalue of the matrix $M$.
A last classical result (see [Shu] for instance) is that : $\log (\lambda( M) ) = h_{\textrm{top}} (\Phi_P)$.
$\square$

 The usefulness of the map  $\Phi_{P}$ is now clear since, thanks to the Markov property, the topological entropy is a computable invariant, as the spectral radius of the Markov transition matrix. 
  In practice the size of the Markov matrix is big but for simple examples the computation is possible.
  
    One advantage of our construction is the fact that the Markov map is well defined for any geometric presentation so it makes the comparison of the volume entropy for different presentations possible.
 
 It turns out that Corollary 5.3 has an immediate consequence:
 
 \begin{lemma}
 If  $P$ is a presentation with $n$ generators (i.e. $| X | = 2 n$ ) and no relations of the form $xxy = Id$, then the following inequalities are satisfied:
 
 {\centerline{  $\log (2n - 3 ) \leqslant h_{top} ( \Phi_{P}) = h_{vol} (P) \leqslant log (2 n -1)  $.}}
 \end{lemma}
 
{\em Proof.} The second inequality is well known, it is just the obvious comparison between the volume entropy of the group presentation with the one for the free group of the same rank. The first inequality is the new result, it is a direct consequence of Corollary 5.3.
 $\square$
 
 The next result is about the special presentations with relations of the form $xxy = Id$, we postpone it's  proof to the Appendix.
 
 \begin{lemma}
 If $P$ is a geometric presentation of $\Gamma$ with $n$ generators and some relations of the form $(*) \hspace{1cm} xxy = Id$,
 then there exists a presentation $P'$ with $n-1$ generators and one relation of the form (*) less than in $P$ so that  $h_{\textrm{vol}}(P') \leqslant h_{\textrm{vol}}(P) $.
 
 \end{lemma}
 
 Finally we obtain :
 
 \begin{theorem}

 The minimal volume entropy of a co-compact hyperbolic surface group is realized, among the geometric presentations, by the presentations with the minimal number of generators.
 \end{theorem}
 
 Lemma 5.14 and 5.15 imply that all the geometric presentations with the minimal number of generators have volume entropy less than the other geometric presentations of the same group. This minimum is realized since the number of such presentations, called {\em minimal}, is finite.
 It remains to prove the following result:
 
 \begin{lemma}
 All the minimal geometric presentations have the same volume entropy.  
 \end{lemma}

 {\em Proof.} Observe that the minimal geometric presentations of co-compact surface groups are very classical, for orientable surfaces of genus $g$ for instance, they have $2 g$ generators and one relation of length $4 g$. There are still several possibilities, for instance in genus $2$ here are two distinct presentations :
 ${\large{< }} a, b, c, d {/} a b a^{-1} b^{-1} c d c^{-1} d^{-1} = id {\large{> } }$  and 
 ${\large{< }} a, b, c, d {/} a b a^{-1}c  d  c^{-1} b^{-1} d^{-1} = id {\large{> } }$.
 
 From the Markov map point of view these two presentations give two different maps but the difference is simply the ordering of the intervals $I_{x_i}$ along the circle $S^1 = \partial \Gamma$.
 All these different maps are constructed with one relation with the same even length and thus all these maps are combinatorially conjugated, in particular the Markov matrices are the same, up to permutation of the indices, 
 and therefore the topological entropies are the same. At the Cayley graph level, the proof is also immediate.$\square$
 
 \subsection{An example.}
 
 For an example, we compute the map $ \Phi_{P}$ for the classical presentation of an orientable surface of genus 2.
  We give some (partial) explicit computations for the geometric presentation : \\
  $ P = < a, b, c, d / a . b. a^{-1} . b^{-1} . c . d . c^{-1} . d^{-1}> $. A small part of the Cayley 2-complex is shown in Figure 10, as well as the subdivision of the interval $I_a$ of the partition. We show bellow the computation for this particular interval and this is enough by the symmetry of the presentation.
 
 \begin{corollary}
 The minimal volume entropy among geometric presentations of genus two surfaces is :
 $log (\frac{3 + \sqrt{17} + \sqrt{ 22 +6 \sqrt{17}}}{2} ) $.
 \end{corollary}

 \begin{figure}[htbp]
\centerline{\includegraphics[height=100mm]{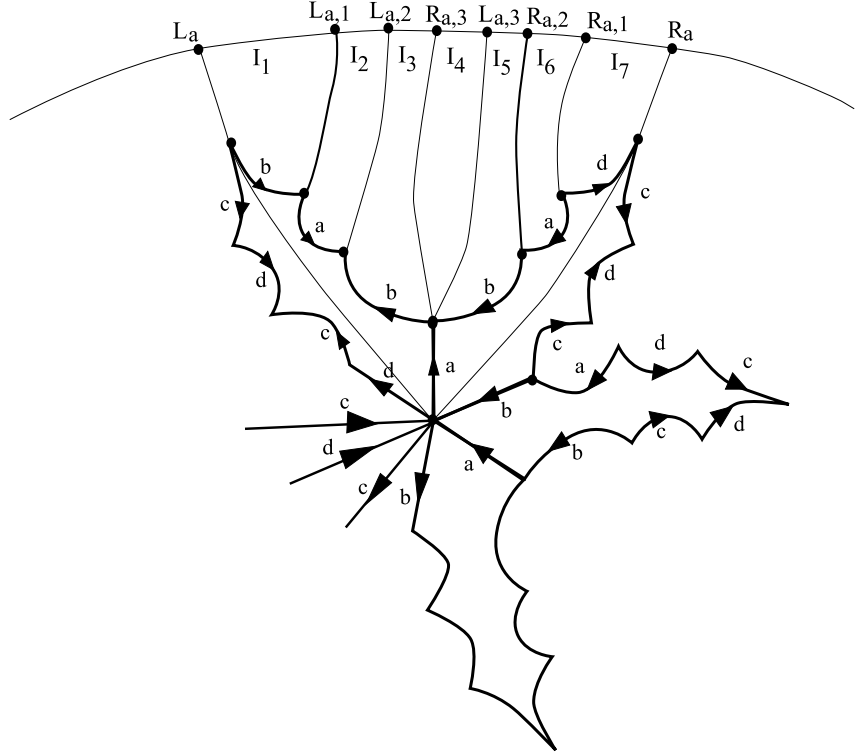} }
\label{fig:10}
\caption{The genus 2 case.}
\end{figure}

 The circle is oriented clockwise and the interval $I_a$ is the concatenation of the 7 intervals 
 $  I_{a, 1} ......  I_{a, 7}$. The subdivision points are denoted :
 $\big\{ L_a, L_{a,1}, L_{a,2}, R_{a,3}, L_{a,3}, R_{a,2}, R_{a,1}, R_{a}\big\}$. With the notations of section 4 these points are written as the limit of the following rays: \\
 
 $ \{L_a \} = a . b. a^{-1} . b^{-1}. \beta^{\infty} [ (b,c)^{opp}], \hspace{1cm} 
  \{R_a \} = a . b^{-1}. a^{-1} . d. \beta^{\infty} [ (c, d^{-1} )^{opp}], \\
  \{L_{a, 1} \} = a . b. a^{-1}.  \beta^{\infty} [ (a, b^{-1})^{opp}],  \hspace{1cm} 
  \{R_{a, 1} \} = a . b^{-1}. a^{-1}. \beta^{\infty} [ (d , a )^{opp}], \\
   \{L_{a, 2} \} = a . b. \beta^{\infty} [ (b^{-1},a^{-1})^{opp}],  \hspace{1.3cm} 
    \{R_{a, 2} \} = a . b^{-1}. \beta^{\infty} [ (a^{-1}, b )^{opp}], \\
    \{L_{a, 3} \} = a .\beta^{\infty} [ (a^{-1}, b )^{opp}],  \hspace{2cm} 
    \{R_{a, 3} \} = a .  \beta^{\infty} [ (b^{-1}, a^{-1} )^{opp}].$
    
    The computation of section 4 gives the following images of these points :
    
$  (1) \hspace {1cm} \Phi_{P} (L_a ) = L_{b, 1}, \hspace {1.2cm}   (5) \hspace {1cm} \Phi_{P} (R_a ) = R_{b^{-1}, 1}, \\
(2) \hspace {1cm} \Phi_{P} (L_{a, 1} ) = L_{b, 2}, \hspace {1cm}  (6) \hspace {1cm} \Phi_{P} (R_{a, 1} ) = R_{b^{-1}, 2}, \\
(3) \hspace {1cm} \Phi_{P} (L_{a, 2} ) = L_{b, 3},  \hspace {1cm}  (7) \hspace {1cm} \Phi_{P} (R_{a, 2} ) = R_{b^{-1}, 3}, \\
(4) \hspace {1cm} \Phi_{P} (L_{a, 3} ) = L_d ,  \hspace {1.2cm} (8) \hspace {1cm} \Phi_{P} (R_{a, 3} ) = R_{d^{-1}}.$ 
 
 The image of the extreme points of an interval gives the image of the interval so we obtain:
 
 $ (1) - (2) \hspace {1cm}  \textrm{gives} :  \hspace {1cm}  \Phi_{P} (I_{a, 1} )  = I_{b, 2}, \\
  (2) - (3) \hspace {1cm} \textrm{gives} :  \hspace {1cm}  \Phi_{P} (I_{a, 2} )  = I_{b, 3}. I_{b, 4}, \\
  (3) - (8) \hspace {1cm} \textrm{gives} :  \hspace {1cm}  \Phi_{P} (I_{a, 3} )  = I_{b, 5}...... I_{d^{-1}, 7}, \\
  (8) - (4) \hspace {1cm} \textrm{gives} :  \hspace {1cm}  \Phi_{P} (I_{a, 4} )  = I_{c^{-1}, 1}...... I_{c^{-1}, 7}, \\
  (4) - (7) \hspace {1cm} \textrm{gives} :  \hspace {1cm}  \Phi_{P} (I_{a, 5} )  = I_{d, 1}...... I_{b^{-1}, 3}, \\
 (7) - (6) \hspace {1cm} \textrm{gives} :  \hspace {1cm}  \Phi_{P} (I_{a, 6} )  = I_{b^{-1}, 4}. I_{b^{-1}, 5}, \\
(6) - (5) \hspace {1cm} \textrm{gives} :  \hspace {1cm}  \Phi_{P} (I_{a, 7} )  = I_{b^{-1}, 6}.$

 The computation of the Markov matrix is now tedious and long, it gives a 56 by 56 integer matrix with many repeating blocks. The computation of the largest eigenvalue has no interest in itself and can be done with any computer. It is interesting to notice that for this classical  presentation, the map defined by Bowen and Series in [BS] is a little bit different and less symmetric than the one presented here, but the number of subdivision intervals is the same (56 in this case) and the computation of the entropy gives the same value. $\square$
 
 The characterisation of the minimal volume entropy of surface groups needs some more steps and some new ideas. It is natural to conjecture that the geometric minimum obtained in this paper is a candidate to be the absolute minimum.\\
 Another class of questions would be to understand more properties of the special circle maps that are defined here. 
 Of course the real challenge would be to define similar maps for arbitrary presentation of an arbitrary hyperbolic group and in particular when the boundary is a higher dimensional sphere.

 \section{ Appendix: Presentations with a relation  xxy = Id.}
 
This Appendix is a proof of  Lemma 5.15. It is treated separately since a new strategy is necessary. The cases where the presentation $P$ has a relation of the form $xxy = Id$ is not covered by Lemmas 5.1 and 5.2 and the conclusion of Corollary 5.3 and Lemma 5.14 are wrong in this case. We assume that $P$ has a relation of the form $(*) \hspace{1cm} xxy = Id$. The construction of the map $\Phi_P$ is the same than in section 4 and we consider the local structure of the 2-complex $Cay^2 (\Gamma, P)$ in this particular case (see Figure 11).
 
  \begin{figure}[htbp]
\centerline{\includegraphics[height=100mm]{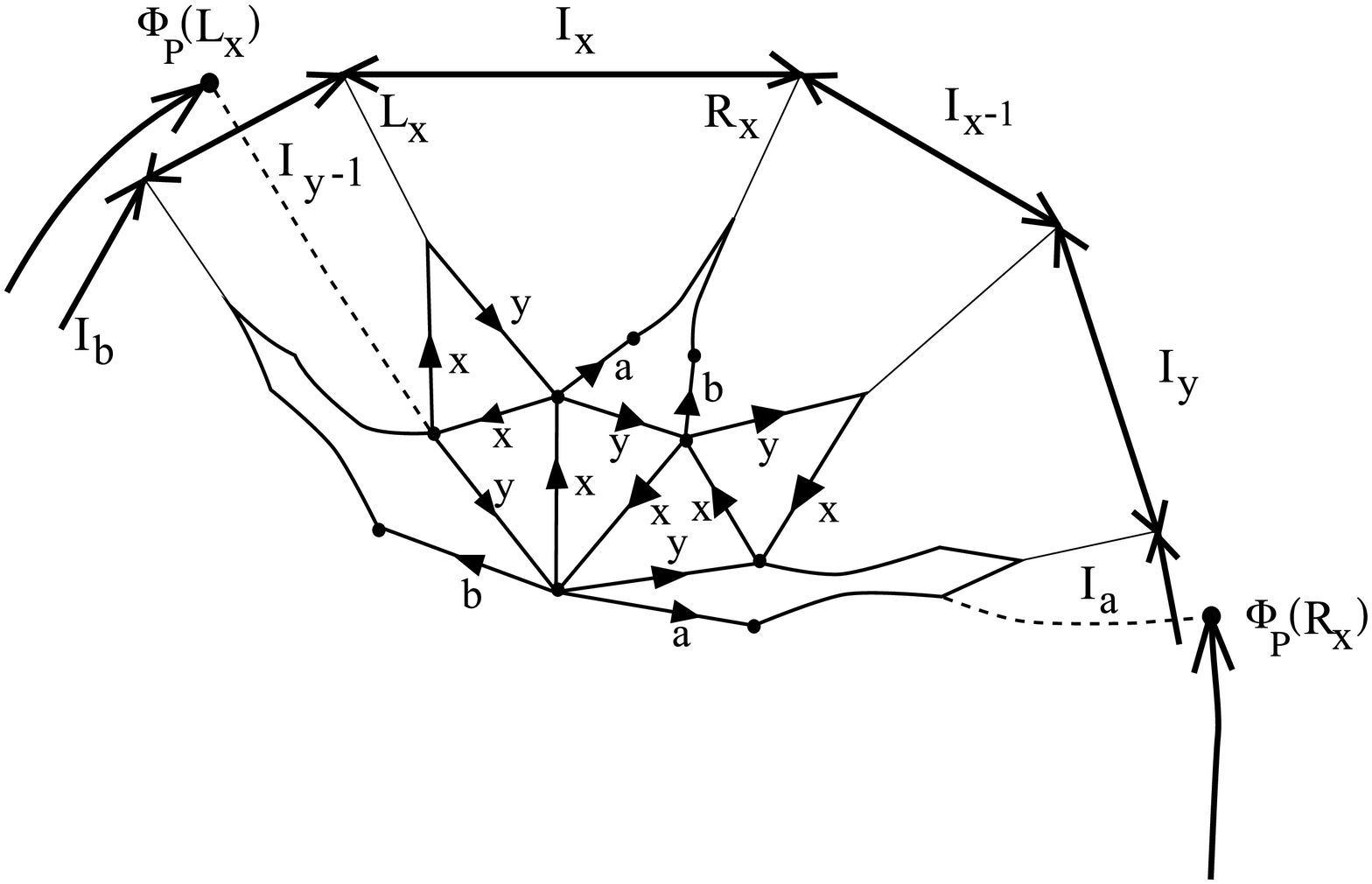} }
\label{fig:11}
\caption{The particular case $xxy = Id$.}
\end{figure}

We consider the partition intervals defining the map $\Phi_P$ with special attention to the 4 intervals 
$I_{x^{\pm 1}},  I_{y^{\pm 1}}$. For the other intervals the conclusions of the Lemmas 5.1 and 5.2 are valid.
We observe in particular that $\Phi_P (I_{x})$ misses 3 intervals  : $I_x$, $I_{x^{-1}}$, $I_y$  as well as a sub-interval before and after (the same is true for $\Phi_P (I_{x^{-1}})$ ). The conclusion of Lemma 5.2 does not apply for these two intervals but it does for  $I_y$ and $I_{y^{-1}}$.

The new step is to transform the presentation $P$ by a specific Dehn twist that preserves the geometric structure, i.e. by a surface Dehn twist (as opposed to a free group Dehn twist). Algebraically this twist is given by  $\tau : P \longrightarrow P '$ such that :\\
\centerline {$\tau ( y) = yx := z$ and $\tau (x_i) = x_i , \forall x_i \in X - \{ y \}$. }

For the new presentation, only two relations are transformed (according to the Figures 11 and 12):\\
\centerline{ $\big\{ xxy, a^{-1}.y.b ...,\ast,\ast, \ast ...\big\} \longrightarrow  \big\{ xz, a^{-1}.z.x^{-1}.b ...,\ast,\ast, \ast  ...\big\} $ .}

We check that $P '$ is geometric since the Dehn twist $\tau$ is realized on the surface and we compute the map $\Phi_{P '}$ (see Figure 12).

\begin{lemma}
With the above notations the topological entropy satisfies :  
$h_{\textrm{top}} ( \Phi_P) \geqslant h_{\textrm{top}} ( \Phi_{P '}) $.

\end{lemma}

{\em Proof.} The number of generators is the same for $P$ and $P '$ and the automorphism $\tau$ gives an identification between the generators that induces an identification of the intervals $I_{x_i}$. 
The partition for the two Markov maps are different, in particular because the relations have different length, but each interval $I_{x_i}$ has a rough partition obtained from the proof of Lemma 5.1 (and 5.2). Indeed, each interval $I_{x_i}$ has a left part $I_{x_i}^L$, a central part $I_{x_i}^C$ and a right part 
$I_{x_i}^R$. These particular partition points are the image under $\Phi_{P }$ of the extreme points : 
$\Phi_{P } (\partial I_{x_i})$, $ x_i \in I_{x_i}$.\\
The proof of Corollary 5.3 is based on the fact that, in the cases covered by Lemma 5.1, each point in these sub-intervals have $2n-1$ pre-images for $I_{x_i}^C$ and $2n-2$ pre-images for $I_{x_i}^L$ and $I_{x_i}^R$. In these cases we say that the corresponding interval $I_{x_i}$ is of {\em type}
$(2n-2, 2n-1, 2n-2)$. For the cases covered by Lemma 5.2 all the intervals $I_{x_i}$ such that $x_i$ does not belong to a relation of length 3 are also of type $(2n-2, 2n-1, 2n-2)$ and the generators $x_i$ that belong to a relation of length 3 define intervals of type $(2n-2, 2n-2, 2n-3)$ (or $(2n-3, 2n-2, 2n-2)$).

For the map $\Phi_{P }$, where $P$ has a relation of the form (*) $xxy = Id$, we observe that the intervals 
$ I_{x_i}, x_i \notin \{ x^{\pm 1} \}$ are of type $(2n-2, 2n-1, 2n-2)$ or $(2n-2, 2n-2, 2n-3)$ and the special intervals $I_{x^{\pm 1}}$ are of type $(2n-3, 2n-3, 2n-4)$. These values are obtained by direct checking (see Figure 11).\\
For the map $\Phi_{P '}$ the situation is quite different. We observe, with the above notations (see Figure 12), that the image $\Phi_{P '} (I_{x^{\pm 1}})$ covers $n-1$ intervals $ I_{x_i}$ minus one sub-interval 
($I^L$ or $I^R$) on one side and $\Phi_{P '} (I_{z^{\pm 1}})$ also covers $n-1$ intervals $ I_{x_i}$ minus one sub-interval ($I^L$ or $I^R$) on one side. Another particular property of this map is that 
$\Phi_{P '} (R_x) = \Phi_{P '} (L_{z^{-1}})$ and  $\Phi_{P '} (L_z) = \Phi_{P '} (R_{x^{-1}})$. 
These observations imply that each interval $ I_{x_i}, x_i \notin \{ x^{\pm 1},  z^{\pm 1}\}$ that was of type $(a, b, c)$ for $\Phi_{P }$ is now of type $(a-2, b-2, c-2)$ for $\Phi_{P '}$. The last 4 intervals 
$I_{x^{\pm 1}}$ and $I_{z^{\pm 1}}$ are of respective types $ ( 2n-4, 2n-3, 2n-3)$ and $ ( 2n-3, 2n-3, 2n-4)$.
This imply that the number of pre-images grow much slower, under iteration by $\Phi_{P '}$ than for 
$\Phi_{P }$. This completes the proof of Lemma 6.1.$\square$

This example shows that an elementary transformations of the presentation, like a Dehn twist here, may have a very strong impact on the entropy even when the number of generators is fixed.

  \begin{figure}[htbp]
\centerline{\includegraphics[height=100mm]{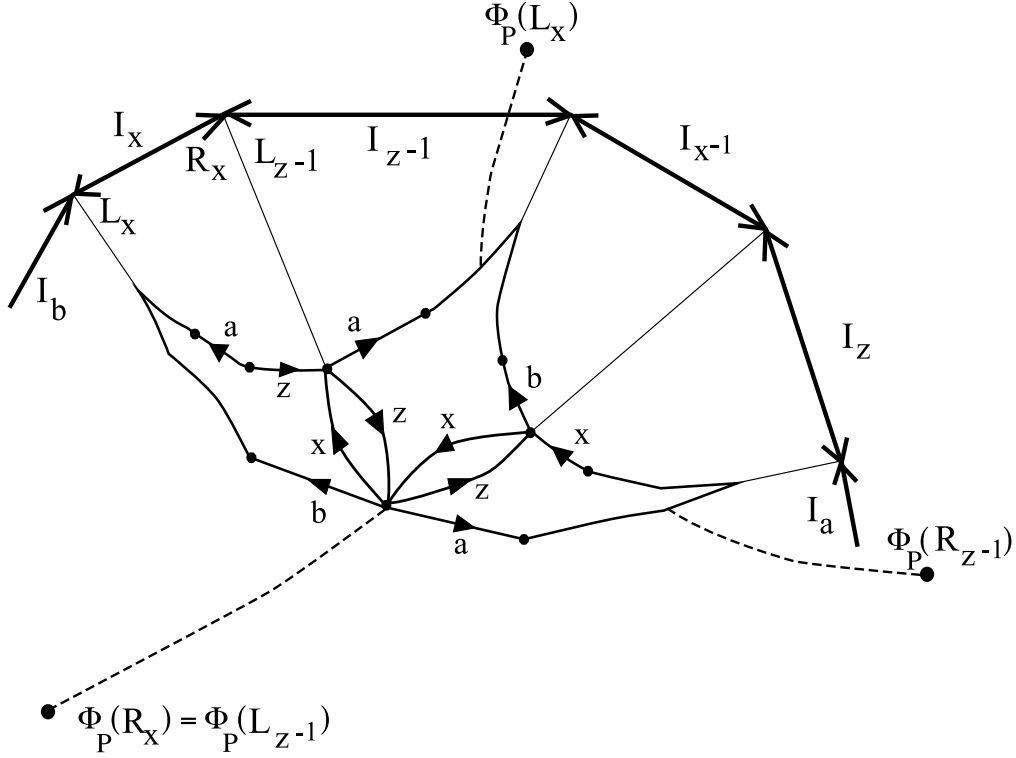} }
\label{fig:12}
\caption{New presentation with a length 2 relation.}
\end{figure}

 The next step is to remove the unnecessary generator ($x$ or $z$).
 
 \begin{lemma}
 If $P '$ is a presentation with $n$ generators and a relation of length 2 :\\
 $(**) \hspace{5mm} x z = Id$, then the new presentation $P ''$ obtained from $P '$ by removing the generator $z$ (for instance) and the relation (**) satisfies : $h_{\textrm{top}} ( \Phi_{P '}) =  h_{\textrm{top}} ( \Phi_{P ''}) $.
 \end{lemma}
 
 {\em Proof.} The bigon $\beta(x, z^{-1}) $ given by the relation of length 2 defines a partition point 
 $(x, z^{-1})^{\infty}$ that has two sides, one is called $R_x$ and the other $L_{ z^{-1}}$ for the two adjacent intervals $I_x$ and $I_{ z^{-1}}$. The main observation is now that :
  $\Phi_{P '}(R_x) = \Phi_{P '} (L_{ z^{-1}})$ (Figure 12 shows an example of this fact). The map 
  $\Phi_{P '}$ is thus continuous at  $(x, z^{-1})^{\infty}$ and
 $\Phi_{P '}$ is a homeomorphism on the larger interval $I_x \bigcup I_{ z^{-1}}$. Therefore we can remove the partition point $(x, z^{-1})^{\infty}$ without changing the dynamics. Removing the partition point corresponds, for the presentation, to replace the 2-cell corresponding to the relation (**) by a single edge, for instance $x$. Since the two dynamics before and after removing the partition point are the same the entropy are the same. $\square$\\
 This completes the proof of Lemma 5.15.

 \pagebreak

 {\bf References.}
 
 [AKM] R.Adler, A. Konheim, M.McAndrew. {\em Topological entropy} Trans. Amer. Math. Soc.  {\bf  114} (1965),  309-319.
 
 [BCG] G.Besson, G.Courtois, S.Gallot. {\em Minimal entropy and Mostow's rigidity theorems.} Ergo. Theo.  Dyn. Syst. 16, (1996), 623-649.
 
 [Bo] R.Bowen. {\em Equilibrium states and the ergodic theory of Anosov diffeomorphisms.} Lec.Notes. Math. Springer, 470, (1975).
 
 [Bou] M. Bourdon. {\em Action quasi-convexes de groupes hyperboliques: flot g\'eod\'esique.} Thesis, Universit\'e Paris Sud XI, Orsay, (1993).
 
 [BS] R.Bowen, C.Series. { \em Markov maps associated with fuschian groups.} Pub.Math. IHES.  {\bf 50}, (1979), 153-170.
 
 [FL] W.J.Floyd, S.P. Plotnick. {\em Growth functions on fuschian groups and the Euler characteristic.}
 Invent.Math. {\bf 88}, (1987), 1-29.
 
 [Gr1] M.Gromov. {\em Hyperbolic groups.} In Essays in group theory, Vol 8 of Math. Sci. Res. Inst. Publ.
 pp 75- 263. Springer, Nw York, (1987).
 
 [Gr2] M.Gromov. {\em Groups of polynomial growth and expanding maps.} Publ.Math. IHES. {\bf 53}
 (1981), 53-78.
 
 [Gri] R.I. Grigorchuk. {\em On Milnor's problem on group growth.} Soviet. Math.Dokl. {\bf 28}, N¡ 1 (1983), 23-26.
 
 [dlH] P.de la Harpe. {\em Topics in geometric group theory.}  Chicago Lectures in Math. (2000).
 
 [He]  G.A.Hedlund. {\em On the metrical transitivity of the geodesics on closed surfaces of constant negative curvature.} Ann.Math. {\bf 35}, (1934), 787-808.
 
 [Mo] G.D. Mostow. {\em Quasi-conformal mappings in n-space and the rigidity of hyperbolic space forms.} Publ. Math. IHES. {\bf 34} (1968), 53-104.
 
 [Ro]  A.C. Rocha. {\em Symbolic dynamics for Kleinian groups.} Thesis, Warwick University, (1994).
 
 [Sho] H.Short et al {\em Notes on word hyperbolic groups.}  In Group theory from a geometric viewpoint.
 Ed E.Ghys, A.Haefliger, A.Verjovsky, World Scientific, (1991), 3-63.
 
 [Shu]  M. Shub. {\em Stabilit\'e globale des syst\`emes dynamiques.} Asterisque, Vol 56, (1978).

\end{document}